\documentclass[final,onefignum,onetabnum]{siamart251216}
 
\usepackage{amssymb}
\usepackage{bm}
\usepackage{physics}
\usepackage{booktabs}
\usepackage{subcaption}
\usepackage[normalem]{ulem}
\usepackage{algpseudocode}   

\DeclareMathOperator*{\argmin}{arg\,min}
\DeclareMathOperator{\diag}{diag}
\DeclareMathOperator{\MAP}{MAP}

\newsiamremark{remark}{Remark}

\headers{Structure-Preserving DA for Conservation Laws}{Y. Xiao and A. Gelb}

\title{Preserving Structure  in Variational Data Assimilation of Hyperbolic
Conservation Laws via Regularization\thanks{\today
\funding{This work was supported by the DOD (ONR MURI) \#N00014-20-1-2595 and DOE ASCR \#DE-SC0025555.}}}

\author{Yao Xiao\thanks{Department of Mathematics, Dartmouth College,
Hanover, NH 03755 (\email{yao.xiao@dartmouth.edu}).}
\and Anne Gelb\thanks{Department of Mathematics, Dartmouth College,
Hanover, NH 03755 (\email{annegelb@math.dartmouth.edu}).}}

\begin{document}
\maketitle
 
\begin{abstract}
Hyperbolic conservation laws pose a challenging setting for variational data assimilation, since the Gaussian assumption imposes a smoothness that smears jump discontinuities. A sparsity promoting regularization can help to mitigate this problem by incorporating a structural prior into the variational objective.  Care must be taken when choosing this prior because it serves to counterbalance both a structurally misspecified background-error covariance, which  introduces spurious oscillations over smooth regions, as well as spatially sparse observations, which leave too little data to constrain the state. Moreover, state variables of hyperbolic conservation laws typically do not have piecewise constant structure, which is an assumption made when using standard  sparsity-promoting operators, such as total variation. Indeed the underlying {\em variability} of hyperbolic conservation law solutions are neither known in advance nor fixed. Hence using higher order total variation is also not suitable. Here we introduce a new regularization term for data assimilation, built on the sparsity promoting residual transform that neither commits to a fixed smoothness order nor requires the underlying variability to be known a priori. We solve the resulting nonconvex objective within a three-dimensional variational framework using generalized sparse Bayesian learning. This approach replaces a global \(\ell _{1}\) penalty with a hierarchical conditional Gaussian prior on the residual transform, where a learned hyper-prior yields location-specific weights attached to each spatial point. Numerical experiments across scalar, shallow-water, and Euler test problems show that this regularization term resolves discontinuities more accurately than total variation, and that its advantage grows as the background-error covariance becomes misspecified and the observations become sparser.
\end{abstract}



\section{Introduction}
Nonlinear hyperbolic conservation laws develop discontinuities in finite time, even from smooth initial data. Shocks and contact discontinuities are therefore intrinsic to the solutions of such systems, and their location and magnitude carry important physical information about the state. The accurate recovery of these features poses the central challenge for data assimilation of conservation laws. Data assimilation methods fall broadly into variational and sequential approaches \cite{evensen2009data,law2015data}. 
Variational methods formulate assimilation as the minimization of a cost function, whereas ensemble-based methods approximate the state distribution using an ensemble of samples. Both rely on a Gaussian assumption for the error statistics, which is well known to be poorly suited to state variables that admit discontinuous profiles. In particular, variational methods penalize the departure of the analysis from a background state through a quadratic term weighted by the background-error covariance. This penalty suppresses sharp variations, so the analysis tends to smear the discontinuities in the state profile. When the background covariance is simplified or poorly conditioned, it further admits spurious oscillations across the domain \cite{freitag2013resolution,ebtehaj2014variational}. These deficiencies are rooted in the underlying Gaussian assumption of a globally smooth state, which is inappropriate for the piecewise smooth solutions of conservation laws. 
To relax the assumption of a globally smooth state, the $\ell_2$ objective can be augmented with a sparsity promoting penalty when the state admits a sparse representation. For example, $\ell_1$ regularization promotes sparsity and is widely used in image recovery, where the image is sparse in the edge domain, and it has since been adopted within the variational data assimilation framework \cite{freitag2010l1, asadi2019data}. When the state variable is sparse in the gradient domain, a total variation penalty on the gradient, or a $\ell_{1,2}$ penalty, resolves jump discontinuities more accurately than Tikhonov regularization does \cite{budd2011regularization, freitag2013resolution}. Such penalties are best suited to piecewise constant states, for which the gradient is genuinely sparse. On the smooth portions of a field a total variation penalty instead produces staircasing artifacts, which motivated the use of total generalized variation \cite{de2019total}. Structure-preserving variational assimilation has been applied across a range of settings, including advection \cite{freitag2013resolution} and shallow-water flows \cite{bernigaud2024non}, hydrological states \cite{ebtehaj2013variational}, and sea ice \cite{asadi2019data}. More broadly, non-quadratic regularization has been used in variational assimilation for related purposes, including $\ell_1$ and Huber norms for robustness to outliers \cite{rao2017robust}, dictionary-learning sparsity for image and flow retrieval \cite{li2021assimilation}, and $\ell_p$-norm penalties \cite{bernigaud2021p}.These methods were developed to improve the recovery of an initial condition over a single assimilation window, using a correctly specified background-error covariance.

These methods have been developed primarily for densely observed states. In practice, however, the number of observations is typically far smaller than the dimension of the state, and their spatial distribution is highly nonuniform \cite{kalnay2003atmospheric,carrassi2018data}. In many applications dense observation is simply not attainable. The subsurface ocean, for instance, is observed only by a sparse array of profiling floats \cite{roemmich2009argo}. The assimilation must therefore reconstruct the state from spatially sparse data. At locations where no data is available, the assimilation is governed by the assumptions imposed on the solution, both through the forecast background and through the structural regularization. The sparser the observations, the more the reconstruction relies on these assumptions. This effect is most pronounced near discontinuities, whose location and jump magnitude must often be inferred where no data is available \cite{li2024structurally}. This sparse regime motivates our use of the variational approach, which encodes the assumptions on the solution explicitly, rather than inferring them from sample statistics. In this work we examine this regime for both single-time and cycled assimilation.

In these methods the background-error covariance is assumed correctly specified, taken to be either diagonal \cite{budd2011regularization, de2019total} or given a prescribed correlation structure \cite{freitag2013resolution,ebtehaj2014variational}. In realistic settings its assumed correlation structure may not match that of the true background error. In this work we explore assimilation under sparse observations and across a range of background-error covariances, from well-specified correlated and diagonal structures to a misspecified one. A similar misspecification has been studied for the observation-error covariance \cite{stewart2013data}. We work within a three-dimensional variational (3D-Var) framework, augmenting a Gaussian background term with a sparsity-promoting prior in the edge domain. For this sparse prior we use the residual prior transform \cite{xiao2026new,xiao2026joint}. Rather than assuming a fixed type of variability, such as piecewise constant, it is formed from the difference of two operators that respond to the state in a similar way, so that the difference is small except for jump discontinuities. We compare it against the standard first order local differencing. The choice of prior is critical. A prior that commits to a fixed smoothness order, such as the piecewise constant assumption of first-order local differencing, fails to capture profiles with substantial variability. Moreover, it requires the smoothness order a priori, which is rarely available in practice. The residual prior transform \cite{xiao2026joint,xiao2026new} is used here to avoid the impact of both limitations on the framework, since it neither assumes a fixed smoothness order nor requires it to be known a priori. 
The resulting objective is solved by generalized sparse Bayesian learning (GSBL), which we adapt from image reconstruction \cite{glaubitz2023generalized} to the data assimilation setting. We evaluate the method on a set of shock-forming conservation laws, including the Burgers, shallow-water, and Euler equations, under both single-time and cycled assimilation.

The rest of this paper is organized as follows. \Cref{sec:prelim} reviews the variational method, the sparse Bayesian learning algorithm, and the residual prior transform. \Cref{sec:sblvda} details how these ingredients are combined to improve the data assimilation. \Cref{sec:B-design} then explains how the background-error covariance is specified for our test cases. \Cref{sec:numerics} evaluates the performance of our method on four hyperbolic conservation law problems. Finally, \cref{sec:conclusion} contains some concluding remarks.

\section{Preliminaries}
\label{sec:prelim}

We consider the state estimation problem for a spatially discretized dynamical system. Let $\bm{x}_t \in \mathbb{R}^n$ denote the system state at discrete time $t$ on a uniform grid of $n$ cells over the spatial domain $[a, b]$, with cell width
\begin{equation}\label{eq:discrete_grid}
    s_j=a+j\Delta s, \quad \Delta s = \frac{b-a}{n}, \quad
    j=0,\dots,n-1,
\end{equation}
so that the $i$th entry $x_i=\bm{x}_t(s_i)$ of $\bm{x}_t$ is the state value at the $i$th cell center. While this investigation focuses on one-dimensional problems for clarity and conciseness, there are no inherent limitations to our approach. Specifically, the residual transform operator~\cite{xiao2026joint} is, by construction,  a pointwise scalar comparison rather than a combination of oriented derivative components.  It therefore does not depend on the orientation of the underlying discontinuity relative to the grid. This avoids the isotropic/anisotropic combination question that arises in multidimensional TV, where the individual derivative components being combined are themselves orientation-dependent. Extending this methodology to higher dimensions is straightforward.  
When necessary to advance the state between observation times, we evolve it with a generally nonlinear forecast operator $\mathcal{M}\colon \mathbb{R}^n \to \mathbb{R}^n$,
\begin{equation}\label{eq:forecast_model}
\bm{x}_{t} = \mathcal{M}(\bm{x}_{t-1}), \quad t = 1, \dots, T,
\end{equation}
where $\bm x_0$ denotes the initial condition of the dynamical system and $t$ indexes the assimilation cycles up to the final cycle $T$. Each application of $\mathcal{M}$ advances the state from one cycle to the next, at which observations are collected. The assimilation is performed within the physical time interval $[0,t_f]$, where $t_f$ denotes the final time, reached at the final cycle $T$. The numerical model may take finer internal time steps within each cycle.

At each assimilation cycle, a noisy observation $\bm{y}_t \in \mathbb{R}^m$ is collected through a linear observation operator $H \in \mathbb{R}^{m \times n}$, where $m \leq n$,
\begin{equation}\label{eq:observation_model}
\bm{y}_t = H \bm{x}_t + \bm{\epsilon}_t, \quad \bm{\epsilon}_t \sim \mathcal{N}(\bm{0}, \Sigma_t),
\end{equation}
where $\Sigma_t$ denotes the observation-error covariance.

Data assimilation combines the observation $\bm y_t$ with the background $\bm x_t^b$, the forecast state estimate before the observation is assimilated, to produce an analysis $\bm x_t^a$. The analysis balances the background against the observation according to their respective uncertainties.
We refer to producing one analysis as an \emph{analysis step}, the basic unit of assimilation, performed once per observation time. In particular, we consider two settings. In the \emph{cycled} setting, an analysis step is carried out at each of the observation times $t=1,\dots,T$ with $T>1$, and the forecast operator $\mathcal M$ propagates the analysis at one time to the background at the next. In this way the background carries forward the information from earlier observations. The \emph{single-cycle} setting is the special case $T=1$, in which the single background is propagated from the initial condition. In both settings the assimilation reaches the same final time $t_f$.

\begin{remark}\label{re:cycle-setup}
When the background is inaccurate, the benefit of assimilation is demonstrated by the single-cycle case, which yields substantial error reduction. On the other hand, the fresh background perturbation injected at every cycle causes the degradation over multiple cycles. This per-cycle injection is the additive inflation used to represent model error in ensemble data assimilation \cite{whitaker2012evaluating}. Unlike multiplicative inflation, which merely rescales the existing spread, it introduces new error structure that models a continually regenerated forecast error. The number of cycles thus measures the ability of each method to suppress this repeated injection.
\end{remark}

We adopt the following notational conventions. For a symmetric positive definite matrix $C$, we denote the weighted norm by $\|\bm{z}\|_{C}^2 = \bm{z}^T C^{-1} \bm{z}$. 

Our approach combines a three-dimensional variational (3D-Var) framework with hierarchical sparse Bayesian learning, using the residual prior transform as the sparsity-promoting operator.\footnote{We note that while the variational framework is well-suited for  the sparse data environment considered in this work, an analogous regularization could also be pursued within a sequential filtering framework, since the regularization enters directly as an objective term. This requires a non-Gaussian treatment beyond a standard covariance-based prior, however.} We review the 3D-Var formulation, the sparse Bayesian learning solver, and the edge-domain operator in the following subsections.

\subsection{Three-dimensional variational data assimilation}
\label{sec:prelim-3dvar}

3D-Var estimates the state at each observation time by combining a background with the observation \cite{lorenc1986analysis, kalnay2003atmospheric}, using a prescribed background-error covariance. The treatment of this covariance is a key departure from sequential filtering, where the covariance is instead propagated by the dynamics, as in the Kalman filter \cite{kalman1960new, evensen1994sequential, evensen2009data}. 

A Gaussian background-error covariance is well justified when the underlying forecast-error distribution remains close to Gaussian, but for hyperbolic conservation laws with shocks, uncertainty in a shock's position induces a multimodal forecast distribution that a linear, unimodal covariance structure cannot represent \cite{zhou2026neural}. 
While inflation and localization can mitigate this loss of structure, they primarily correct the scale and correlation shape of an otherwise unimodal Gaussian covariance. They cannot resolve the more fundamental mismatch between a unimodal model and the multimodal uncertainty actually induced by a shock's position. This mismatch becomes unavoidable when observations are sparse.  Specifically, because we cannot rely on nearby data, the analysis is stuck with the limitations of the background error covariance. This motivates our use of a prescribed background-error covariance in 3D-Var, together with explicit regularization, which allows structural priors to be imposed directly, rather than through the  (inadequate) Gaussian covariance alone. 
We now describe its background and analysis steps.

\subsubsection{Background}\label{sec:prelim-background}

The background $\bm x_t^b$ is the forecast state estimate at each cycle $t=1,\dots, T$, obtained by advancing the previous analysis through the forecast model $\mathcal M$ in \eqref{eq:forecast_model} with additive noise,
\begin{equation}\label{eq:background_gen}
\bm{x}_t^b = \mathcal{M}(\bm{x}_{t-1}^a) + \bm{\eta}_t, \quad \bm{\eta}_t \sim \mathcal{N}(\bm{0}, B_t),
\end{equation}
where $\bm{x}_{t-1}^a$ is the analysis state from the previous cycle and $\bm\eta_t$ is zero-mean Gaussian noise with background-error covariance $B_t\in\mathbb R^{n\times n}$.
In the single-cycle setting, the background is formed from the forecast $\mathcal M(\bm x_0)$ propagated from the initial time to the analysis time with additive noise. The covariance $B_t$ specifies the assumed statistics of this background error and acts as the weighting matrix for the background term in the analysis objective. Its structure is the core experimental design choice in this investigation, and we defer its specification to \cref{sec:B-design}.

\subsubsection{Analysis}\label{sec:prelim-analysis}

At each cycle $t=1,\dots,T$, the analysis is obtained by minimizing the 3D-Var cost function \cite{lorenc1986analysis, courtier1994strategy, kalnay2003atmospheric, carrassi2018data},
\begin{equation}\label{eq:3dvar_cost}
\mathcal G(\bm{x}_t) = \frac{1}{2}\norm{\bm{x}_t - \bm{x}_t^b}_{B_t}^2 + \frac{1}{2}\norm{\bm{y}_t - H\bm{x}_t}_{\Sigma_t}^2,
\end{equation}
a least-squares measure of the distance of the state from the background and from the observation, weighted by the background-error covariance and observation-error covariance, respectively. The first term penalizes the departure of the state from the background $\bm x_t^b$, whereas the second penalizes the misfit to the observation $\bm y_t$. The quadratic form corresponds to the assumption of Gaussian background and observation errors. 

The analysis state is the minimizer of \eqref{eq:3dvar_cost},
\begin{equation}\label{eq:3dvar_analysis}
\bm{x}_t^a = \argmin_{\bm{x}_t} \mathcal{G}(\bm{x}_t).
\end{equation}
For a linear observation operator $H\in\mathbb R^{m\times n}$, the cost function is quadratic and the analysis solves the normal equations
\begin{equation}\label{eq:3dvar_normal}
\left( B_t^{-1} + H^T \Sigma_t^{-1} H \right) \bm{x}_t^a = B_t^{-1} \bm{x}_t^b + H^T \Sigma_t^{-1} \bm{y}_t,
\end{equation}
where $\cdot^T$ denotes the matrix transpose. The solution can equivalently be written as the Kalman filter analysis with gain $B_t H^T (H B_t H^T + \Sigma_t)^{-1}$, which requires inverting only the $m\times m$ matrix $H B_t H^T + \Sigma_t$ in the observation space, rather than an $n\times n$ system in the state space \cite{lorenc1986analysis, kalnay2003atmospheric}.

\subsection{Regularized 3D-Var}
\label{sec:prelim-RVDA}

Physical variables often exhibit piecewise smooth behavior in space. The quadratic cost function \eqref{eq:3dvar_cost} promotes smooth variations in the state and smears the jump discontinuities that form under a conservation law. Moreover, when the observation is severely corrupted by noise, the $\ell_2$ penalty produces highly oscillatory analyses. A piecewise smooth state is sparse in a suitable transform domain, such as the edge domain. To exploit this sparsity and mitigate these issues, a sparsity-promoting penalty in the form of $\ell_1$ regularization has been introduced into the cost function \cite{freitag2010l1, budd2011regularization},
\begin{equation}\label{eq:rvda_cost}
\mathcal{G}(\bm{x}_t) = \frac{1}{2}\norm{\bm{x}_t - \bm{x}_t^b}_{B_t}^2 + \frac{1}{2}\norm{\bm{y}_t - H\bm{x}_t}_{\Sigma_t}^2 + \lambda \norm{\Phi \bm{x}_t}_1,
\end{equation}
where $\Phi\in\mathbb R^{r\times n}$ is a sparsifying operator that maps the state to its representation in the edge domain, and $\lambda>0$ controls the strength of the penalty. The standard total variation operator is a popular choice, favoring piecewise constant structure, and we benchmark the regularized variational results against it. However, as noted in the introduction, physical quantities rarely conform to a single order of smoothness, which is moreover not known in advance. Different smooth regions of a state may vary at different orders, and applying the first order total variation operator to such a highly variable quantity introduces ``staircasing'' artifacts. To circumvent this limitation, we instead use the residual prior transform \cite{xiao2026new,xiao2026joint},  which does not assume a fixed or known order of smoothness. We describe its construction in \cref{sec:prelim_sparseOps}.

The $\ell_1$ regularization is a canonical means of promoting sparsity in sparse data settings, and our residual prior transform is well suited to this framework. Unlike the quadratic cost function \eqref{eq:3dvar_cost}, however, \eqref{eq:rvda_cost} has no closed-form solution and must be minimized iteratively \cite{freitag2010l1, budd2011regularization, freitag2013resolution, ebtehaj2014variational}, owing to the nondifferentiability of the $\ell_1$ norm. The LASSO is a widely used remedy, but it requires a hand-tuned regularization parameter and inherits this same nondifferentiability. We instead adopt generalized sparse Bayesian learning (GSBL) \cite{glaubitz2023generalized} which promotes sparsity through a sequence of differentiable, reweighted subproblems. This gain comes at the cost of a nonconvex objective, which, unlike the convex LASSO, may admit local minima \cite{wipf2004sparse, glaubitz2023generalized}.

\subsection{Hierarchical sparse Bayesian learning}
\label{sub:prelim_gsbl}

The generalized sparse Bayesian learning (GSBL) framework \cite{glaubitz2023generalized} is a hierarchical Bayesian approach \cite{tipping2001sparse, calvetti2020sparse} for recovering a piecewise smooth signal from noisy and undersampled data. It computes a maximum a posteriori (MAP) point estimate of the signal under a sparsity-promoting prior, whose strength is governed by hyperparameters inferred from the data rather than fixed in advance. Moreover, its hierarchical structure yields a posterior covariance that quantifies the confidence in the recovered estimate. These properties make GSBL well suited to signals that combine sharp features with smoothly varying regions, such as the solutions of hyperbolic conservation laws considered here. We give a brief review below to keep the presentation self-contained.

\subsubsection{Hierarchical sparse Bayesian model}
\label{sub:prelim_gsbl_model}

Consider a piecewise smooth function $f \colon [a, b] \to \mathbb{R}$ discretized at $n$ uniformly spaced grid points  in \eqref{eq:discrete_grid}.
We seek to recover the discretized signal $\bm{f} = [f(s_0), \dots, f(s_{n-1})]^T \in \mathbb{R}^n$ from a noisy observation $\bm{d} \in \mathbb{R}^m$ through the forward model
\begin{equation}\label{eq: gsbl_forward}
    \bm{d} = F \bm{f} + \bm{\eta},
\end{equation}
where $F \in \mathbb{R}^{m \times n}$ is a known forward operator and $\bm{\eta}$ is additive, independent and identically distributed (i.i.d.) Gaussian noise with zero mean and precision $\alpha > 0$,
\begin{equation}\label{eq: gsbl_noise}
    \bm{\eta} \sim \mathcal{N}(\bm{0},\, \alpha^{-1} I_m).
\end{equation}
A piecewise smooth signal is assumed to admit a sparse representation in some transformed domain, such as the edge domain, where its nonzero values are concentrated at the discontinuities.

The GSBL framework formulates this recovery through the three-level hierarchical
model
\begin{subequations}\label{eq: gsbl_model}
\begin{align}
    \bm{d} \mid \bm{f} 
        &\sim \mathcal{N}\!\left(F \bm{f},\, \alpha^{-1} I_m\right), \label{eq: gsbl_likelihood}\\
    \Phi \bm{f} \mid \bm{\theta} 
        &\sim \mathcal{N}\!\left(\bm{0},\, D_{\bm{\theta}}^{-1}\right), \label{eq: gsbl_prior}\\
    \theta_k 
        &\sim \Gamma(\beta, \vartheta), \quad k = 1, \dots, K. \label{eq: gsbl_hyperprior}
\end{align}
\end{subequations}
The likelihood \eqref{eq: gsbl_likelihood} follows directly from the forward model \eqref{eq: gsbl_forward}. The prior operator $\Phi \in \mathbb{R}^{K \times n}$ formalizes the sparsity assumption introduced above. It maps the signal into the transformed domain, so that $\Phi \bm{f}$ encodes the structural belief that the signal has few nonzero values away from its discontinuities. The conditionally Gaussian prior \eqref{eq: gsbl_prior} imposes this sparsity through the diagonal precision matrix $D_{\bm{\theta}} = \diag(\bm{\theta}) \in \mathbb{R}^{K \times K}$, whose entries $\bm{\theta} = [\theta_1, \dots, \theta_K]^T$ act as individual precisions for the components of $\Phi \bm{f}$. A large $\theta_k$ forces the $k$-th component toward zero, while a small $\theta_k$ permits a nonzero entry, so that the support of $\Phi \bm{f}$ corresponds to the components assigned a small precision. The hyperprior \eqref{eq: gsbl_hyperprior} models each precision $\theta_k$ as gamma distributed with shape parameter $\beta$ and scale parameter $\vartheta$, governing how strongly sparsity is promoted. The choice $\beta = 1$ and $\vartheta \to \infty$ yields an uninformative hyperprior, under which the prior precisions $\bm{\theta}$ are inferred from the data rather than fixed in advance \cite{glaubitz2023generalized, calvetti2020sparse}.

The densities corresponding to the hierarchical model \eqref{eq: gsbl_model} are
\begin{subequations}\label{eq: gsbl_densities}
\begin{align}
    \pi(\bm{d} \mid \bm{f}) 
        &\propto \exp\!\left( -\frac{\alpha}{2} \norm{ F \bm{f} - \bm{d} }_2^2 \right), 
        \label{eq: gsbl_density_likelihood}\\
    \pi(\bm{f} \mid \bm{\theta}) 
        &\propto \det(D_{\bm{\theta}})^{1/2} 
            \exp\!\left( -\frac{1}{2} \norm{ D_{\bm{\theta}}^{1/2} \Phi \bm{f} }_2^2 \right), 
        \label{eq: gsbl_density_prior}\\
    \pi(\bm{\theta}) 
        &= \prod_{k=1}^{K} \Gamma(\theta_k \mid \beta, \vartheta) 
        \propto \det(D_{\bm{\theta}})^{\beta - 1} 
            \exp\!\left( -\sum_{k=1}^{K} \frac{\theta_k}{\vartheta} \right). 
        \label{eq: gsbl_density_hyperprior}
\end{align}
\end{subequations}
By Bayes' theorem, the joint posterior of the signal $\bm{f}$ and the hyperparameter
vector $\bm{\theta}$ given the observation $\bm{d}$ factors as the product of the
likelihood, prior, and hyperprior,
\begin{equation}\label{eq: gsbl_bayes}
    \pi(\bm{f}, \bm{\theta} \mid \bm{d}) 
        \propto \pi(\bm{d} \mid \bm{f})\, \pi(\bm{f} \mid \bm{\theta})\, \pi(\bm{\theta}).
\end{equation}
Substituting \eqref{eq: gsbl_densities} into \eqref{eq: gsbl_bayes} yields the joint
posterior density
\begin{equation}\label{eq: gsbl_posterior}
    \pi(\bm{f}, \bm{\theta} \mid \bm{d}) 
        \propto \det(D_{\bm{\theta}})^{\beta - 1/2} 
            \exp\!\left( 
                -\frac{\alpha}{2} \norm{ F \bm{f} - \bm{d} }_2^2 
                -\frac{1}{2} \norm{ D_{\bm{\theta}}^{1/2} \Phi \bm{f} }_2^2 
                -\sum_{k=1}^{K} \frac{\theta_k}{\vartheta} 
            \right).
\end{equation}

\subsubsection{MAP estimation via block coordinate descent}
\label{sub:prelim_gsbl_map}

We estimate the signal and hyperparameters jointly through the MAP estimate
$(\bm{f}^{\MAP}, \bm{\theta}^{\MAP})$, defined as the maximizer of the joint posterior
\eqref{eq: gsbl_posterior}, or equivalently the minimizer of its negative logarithm,
\begin{equation}\label{eq: gsbl_map}
    (\bm{f}^{\MAP}, \bm{\theta}^{\MAP}) 
        = \argmin_{\bm{f}, \bm{\theta}} \mathcal{G}(\bm{f}, \bm{\theta}).
\end{equation}
Discarding the terms independent of $\bm{f}$ and $\bm{\theta}$ gives the objective
\begin{equation}\label{eq: gsbl_objective}
    \mathcal{G}(\bm{f}, \bm{\theta}) 
        = \frac{\alpha}{2} \norm{ F \bm{f} - \bm{d} }_2^2 
        + \frac{1}{2} \norm{ D_{\bm{\theta}}^{1/2} \Phi \bm{f} }_2^2 
        + \sum_{k=1}^{K} \frac{\theta_k}{\vartheta} 
        - \left(\beta - \tfrac{1}{2}\right) \sum_{k=1}^{K} \log \theta_k .
\end{equation}

The objective $\mathcal{G}$ is not jointly convex in $(\bm{f}, \bm{\theta})$, but it is convex in each block when the other is held fixed. The GSBL algorithm exploits this structure through block coordinate descent \cite{wright2015coordinate, beck2017first}, alternating between a minimization over $\bm{f}$ for fixed $\bm{\theta}$ and a minimization over $\bm{\theta}$ for fixed $\bm{f}$. Both subproblems are solved efficiently, as summarized in \cref{alg: gsbl}.

\begin{algorithm}[t]
\caption{GSBL via block coordinate descent}
\label{alg: gsbl}
\begin{algorithmic}[1]
\Require forward operator $F$, data $\bm{d}$, noise precision $\alpha$, prior operator
    $\Phi$, hyperprior parameters $\beta, \vartheta$
\Ensure MAP estimate $\bm{f}^{\MAP}$, $\bm{\theta}^{\MAP}$
\State initialize $\bm{\theta} \gets \bm{1}$, $\bm{f} \gets \bm{0}$
\Repeat
    \State $\bm{f} \gets$ solution of
        $\left(\alpha F^T F + \Phi^T D_{\bm{\theta}} \Phi\right) \bm{f} = \alpha F^T \bm{d}$
    \State $\theta_k \gets \dfrac{\beta - 1/2}{[\Phi \bm{f}]_k^2 / 2 + \vartheta^{-1}}$,
        \quad $k = 1, \dots, K$
\Until{convergence}
\State \Return $\bm{f}^{\MAP} \gets \bm{f}$, $\bm{\theta}^{\MAP} \gets \bm{\theta}$
\end{algorithmic}
\end{algorithm}

For fixed $\bm{\theta}$, minimizing \eqref{eq: gsbl_objective} over $\bm{f}$ reduces to
the quadratic problem
\begin{equation}\label{eq: gsbl_f_update}
    \bm{f}^{\MAP} 
        = \argmin_{\bm{f}} \left\{ 
            \alpha \norm{ F \bm{f} - \bm{d} }_2^2 
            + \norm{ D_{\bm{\theta}}^{1/2} \Phi \bm{f} }_2^2 
        \right\},
\end{equation}
whose minimizer solves the normal equations
\begin{equation}\label{eq: gsbl_normal}
    \left( \alpha F^T F + \Phi^T D_{\bm{\theta}} \Phi \right) \bm{f} = \alpha F^T \bm{d}.
\end{equation}
The solution is unique provided $\ker(F) \cap \ker(\Phi) = \{\bm{0}\}$, a standard condition in the analysis of ill-posed inverse problems \cite{kaipio2005statistical}. Under this condition the system matrix is symmetric positive definite, and \eqref{eq: gsbl_normal} is solved by Cholesky factorization, or by the preconditioned conjugate gradient method \cite{saad2003iterative} for large problems. 

For fixed $\bm{f}$, the objective \eqref{eq: gsbl_objective} decouples across the
components of $\bm{\theta}$, and each precision is updated independently by minimizing
\begin{equation}\label{eq: gsbl_theta_min}
    \theta_k^{\MAP} 
        = \argmin_{\theta_k > 0} \left\{ 
            \frac{\theta_k [\Phi \bm{f}]_k^2}{2} 
            + \frac{\theta_k}{\vartheta} 
            - \left(\beta - \tfrac{1}{2}\right) \log \theta_k 
        \right\}.
\end{equation}
Setting the derivative to zero yields the closed-form update
\begin{equation}\label{eq: gsbl_theta_update}
    \theta_k^{\MAP} 
        = \frac{\beta - 1/2}{[\Phi \bm{f}]_k^2 / 2 + \vartheta^{-1}}.
\end{equation}
The update reflects the sparsity mechanism of the prior. At a discontinuity $[\Phi \bm{f}]_k$ is large, which drives $\theta_k$ small and assigns low precision to that component, permitting a nonzero entry. In smooth regions $[\Phi \bm{f}]_k$ is close to zero, which drives $\theta_k$ large and forces the corresponding entry toward zero. Following \cite{glaubitz2023generalized}, in our experiments we set $\beta=1$ and $\vartheta=10^{-4}$. 

Conditioned on the converged hyperparameters $\bm{\theta}^{\MAP}$, the posterior of
$\bm{f}$ is Gaussian,
\begin{equation}\label{eq: gsbl_cond_posterior}
    \bm{f} \mid \bm{\theta}^{\MAP}, \bm{d} 
        \sim \mathcal{N}\!\left( \bm{f}^{\MAP},\, \Gamma \right), 
        \qquad 
        \Gamma = \left( \alpha F^T F + \Phi^T D_{\bm{\theta}^{\MAP}} \Phi \right)^{-1}.
\end{equation}
The posterior mean coincides with the MAP estimate $\bm{f}^{\MAP}$ from
\eqref{eq: gsbl_normal}, and the covariance $\Gamma$ quantifies the confidence in the
recovered signal.

\subsection{The sparsity promoting prior}
\label{sec:prelim_sparseOps}
{\em Traditional} priors serve to promote sparsity in some transform domain. The standard first-order total variation (TV) is a perfect fit for piecewise constant functions, while piecewise linear or more variable functions require priors such as higher-order total variation (HOTV). These priors commit to a {\em fixed} smoothness order and assume it is known a priori. Neither holds for solutions of hyperbolic conservation laws, whose smoothness varies across regions and whose order is rarely known in advance. While raising the order accommodates more variability, it also introduces oscillations of its own. The residual prior transform \cite{xiao2026new,xiao2026joint} removes the assumption altogether, adapting to variable smoothness without a prescribed order. For self-containment purposes, we now provide a brief description.

The residual prior transform $n\times n$ matrix can be generally defined as 
\begin{equation}\label{eq:residual}
R = L_1 - L_2, 
\end{equation}
where $L_1{\bm f} \approx L_2{\bm f}$ in some metric but $L_1 \ne L_2$.  We note that neither operator used to construct $R$ needs to be sparsity promoting. Any pair of operators that provide a pointwise scalar comparison can be used provided they agree over the smooth regions and respond similarly at the discontinuities of the underlying signal. 

While there are a variety of ways to choose $L_1$ and $L_2$ (see \cite{xiao2026new,xiao2026joint}), here we conveniently prescribe them as edge detectors, and more explicitly as the TV operator \cite{rudin1992nonlinear} and a first-order (pseudo-)spectral Fourier edge detector \cite{gelb2002spectral}, respectively.  In this regard, we define the corresponding edge function to piecewise smooth $f$ as
\[[f](s) = f(s^+) - f(s^-),\] 
which vanishes wherever $f$ is smooth. We then define the ground truth edge vector $\bm{g} \in \mathbb{R}^n$ by $g_j = [f](s_j)$ at the grid points \eqref{eq:discrete_grid}. This vector is sparse, with nonzero entries only at the cells containing a discontinuity. Both  $L_1$ and  $L_2$, are first order edge detectors such that $L_1\bm x\approx L_2\bm x \approx \bm g$.

For purposes of self containment, we explicitly provide
\begin{equation}
    \label{eq:L_1}
L_1 = \begin{pmatrix}
-1 & 1 & 0 & \cdots & 0 & 0 \\
0 & -1 & 1 & \cdots & 0 & 0 \\
0 & 0 & -1 & \cdots & 0 & 0 \\
\vdots & \vdots & \vdots & \ddots & \vdots & \vdots \\
0 & 0 & 0 & \cdots & -1 & 1 \\
1 & 0 & 0 & \cdots & 0 & -1
\end{pmatrix} \end{equation}

To construct $L_2$, we first note that the first-order polynomial pseudo-spectral concentration factor method is equivalent to the centered difference of the Dirichlet kernel \cite{gelb2002spectral}, which we write as 
\begin{equation}
      D(\varphi) := \sum_{k=1}^{n/2} \cos\!\Big(\frac{2\pi k}{n}\varphi\Big)
  = \frac{\sin\!\big(\tfrac{\pi}{2}\varphi\big)}{\sin\!\big(\tfrac{\pi}{n}\varphi\big)}\,
    \cos\!\Big(\frac{\pi}{2}\varphi + \frac{\pi}{n}\varphi\Big).
\end{equation}
The entries of $L_2$ then take the form 
\begin{equation}\label{eq:g_tau}
  g_\tau(m) = \frac{2}{n}\Big[D\big(\mu+\tfrac12\big) - D\big(\mu-\tfrac12\big)\Big],
  \qquad \mu = m + \tau, \qquad m = j-j'.
\end{equation}
Here $\tau \in (0,\frac{1}{2})$ is an offset parameter allowing $L_2$ to compute the approximation on ${s}_{j+\tau}$.  The resulting matrix $L_2$ is the circulant matrix generated by $g_\tau$:
\begin{equation}
  L_2 = \operatorname{circ}\big(g_\tau(0),\, g_\tau(1),\, \ldots,\, g_\tau(n-1)\big),
  \qquad (L_2)_{j,j'} = g_\tau(j-j'),
\end{equation}
explicitly,
\begin{equation}
  L_2 =
  \begin{pmatrix}
    g_\tau(0)      & g_\tau(-1)     & g_\tau(-2)     & \cdots & g_\tau(-(n-1)) \\
    g_\tau(1)      & g_\tau(0)      & g_\tau(-1)     & \cdots & g_\tau(-(n-2)) \\
    g_\tau(2)      & g_\tau(1)      & g_\tau(0)      & \cdots & g_\tau(-(n-3)) \\
    \vdots         & \vdots         & \vdots         & \ddots & \vdots \\
    g_\tau(n-1)    & g_\tau(n-2)    & g_\tau(n-3)    & \cdots & g_\tau(0)
  \end{pmatrix}.
\end{equation}

\begin{remark}[Offset parameter $\tau$] \label{re:offset}
The offset $\tau$ is chosen to prevent the residual transform $R := L_1 - L_2$
from becoming low rank. If $R$ were low rank, the residual operator would fail
to regularize across most of the spectrum, since it would only act nontrivially
on a small subspace (see \cite{xiao2026new} for a detailed discussion of this
failure mode). In our experiments we choose $\tau = \tfrac14$.
\end{remark}

\section{Regularized variational data assimilation via sparse Bayesian learning}
\label{sec:sblvda}

With all of the necessary ingredients in hand, we are now able to formulate the assimilation procedure carried out at each cycle, combining the 3D-Var objective, the residual prior transform, and the GSBL solver into a single regularized variational computation via sparse Bayesian learning.

In order to solve the assimilation problem with the GSBL algorithm of \cref{alg: gsbl}, we interpret the sparsity-promoting cost function \eqref{eq:rvda_cost} in Bayesian terms. Treating the state as a random vector, its negative logarithm is the negative log-posterior of the state given the observation, so the analysis state coincides with the maximizer of the posterior $p(\bm{x}_t \mid \bm{y}_t) \propto p(\bm{y}_t \mid \bm{x}_t)\, p(\bm{x}_t)$ which is analogous to \eqref{eq: gsbl_bayes}.
The observation term defines the Gaussian likelihood \eqref{eq: gsbl_likelihood} with the noise model \eqref{eq: gsbl_noise}, and the $\ell_1$ penalty defines the sparsity-promoting prior \eqref{eq: gsbl_prior} on $\Phi\bm x_t$. The background term is a Gaussian prior on the state with covariance $B_t$. Note that the background prior has no counterpart in the GSBL model \eqref{eq: gsbl_model}. To solve \eqref{eq:rvda_cost} within the GSBL framework, we fold the two $\ell_2$ terms into a single fidelity, as described next.

We first whiten both observation and background residuals using the Cholesky factors $\Sigma_t^{-1/2}$ and $B_t^{-1/2}$, respectively, giving $\Sigma_t^{-1/2}(\bm{y}_t - H\bm{x}_t)$ and $B_t^{-1/2}(\bm{x}_t - \bm{x}_t^b)$.
Stacking the whitened observation and background residuals gives the augmented forward operator and data
\begin{equation}\label{eq:whitened_stack}
\tilde{F}_t = \begin{bmatrix} \Sigma_t^{-1/2} H \\ B_t^{-1/2} \end{bmatrix}, \qquad
\tilde{\bm{d}}_t = \begin{bmatrix} \Sigma_t^{-1/2} \bm{y}_t \\
B_t^{-1/2} \bm{x}_t^b \end{bmatrix}.
\end{equation}
Since whitening a residual by an inverse square root of its covariance turns its Euclidean norm into the corresponding weighted norm, we obtain
\begin{equation}\label{eq:whitened_fidelity}
\norm{\tilde{F}_t \bm{x}_t - \tilde{\bm{d}}_t}_2^2 = \norm{\bm{y}_t - H\bm{x}_t}_{\Sigma_t}^2 + \norm{\bm{x}_t - \bm{x}_t^b}_{B_t}^2.
\end{equation}
The two quadratic terms of \eqref{eq:rvda_cost} are thereby fused into a single fidelity $\norm{\tilde F_t\bm x_t-\tilde{\bm d}_t}_2^2$, with $B_t$ and $\Sigma_t$ absorbed into $\tilde F_t$ and $\tilde{\bm d}_t$.
This matches the GSBL forward model \eqref{eq: gsbl_forward} with $F = \tilde F_t$, $ \bm d=\tilde {\bm d}_t$ and the signal $\bm f$ identified with the state $\bm x_t$. The whitening leaves unit noise precision $\alpha=1$ in \eqref{eq: gsbl_noise}.

After formulating \eqref{eq:rvda_cost} in the form of one $\ell_2$ fidelity term with one $\ell_1$ regularization term, we apply GSBL to obtain the analysis state at cycle $t$. The GSBL $\bm f$-update \eqref{eq: gsbl_f_update} reads
\begin{equation}
    \label{eq:whiten_cost}
    \bm x_t^a = \argmin_{\bm x_t}\norm{\tilde F_t\bm x_t-\tilde{\bm d}_t}_2^2+\norm{D_{\boldsymbol{\theta}}^{1/2}\Phi\bm x_t}_2^2,
\end{equation}
a quadratic formulation that admits a closed-form solution. Substituting $\tilde F_t$ and $\tilde{\bm d}_t$ into \eqref{eq: gsbl_normal}, the normal equations become
\begin{equation}\label{eq:da_normal}
\left( B_t^{-1} + H^T \Sigma_t^{-1} H + \Phi^T D_{\bm{\theta}} \Phi \right) \bm{x}_t = B_t^{-1} \bm{x}_t^b + H^T \Sigma_t^{-1} \bm{y}_t. 
\end{equation}
Note that \eqref{eq:da_normal} differs from \eqref{eq:3dvar_normal} in the additional term $\Phi^T D_{\bm{\theta}} \Phi$, which imposes a structural prior on the state to promote sparsity in the edge domain. The update of hyperparameter $\boldsymbol{\theta}$ follows \eqref{eq: gsbl_theta_update}.
$\Phi$ is the sparsity-promoting transform, for which we use the residual prior transform $\Phi=R$ in \eqref{eq:residual} and benchmark it against the standard TV prior transform $\Phi = L_1$ in \eqref{eq:L_1}.

The analysis $\bm x_t^a$ is obtained by running the GSBL algorithm \Cref{alg: gsbl} with the augmented forward model $\tilde F_t$, data $\tilde{\bm d}_t$, and prior operator $\Phi$, alternating the $\bm x_t$-update \eqref{eq:da_normal} and the $\boldsymbol{\theta}$-update \eqref{eq: gsbl_theta_update} until convergence. This completes the assimilation at each cycle.

\section{Background-error covariance specification}
\label{sec:B-design}

The background-error covariance $B_t$ plays a dual role in our 3D-Var assimilation framework. Its first role is apparent in \eqref{eq:background_gen}, where it represents the covariance of the background error, that is, the departure of the background from the true PDE profile. It also serves a secondary purpose as the weighting of the background term in the analysis objective \eqref{eq:rvda_cost}, encoding the assumed statistics of that error. 

In our specification we begin with the self-consistent case, in which the two roles coincide. That is, the covariance assumed in \eqref{eq:rvda_cost} equals the true covariance used in \eqref{eq:background_gen}. This is the standard assumption underlying the 3D-Var analysis of \Cref{sec:prelim-3dvar}. 
In real applications, however, the true covariance is rarely known and requires estimation. This leads to a misspecification between the covariance used in \eqref{eq:background_gen} and its counterpart in \eqref{eq:rvda_cost}. Specifically, we denote the true generating covariance by $B_t$ as in \eqref{eq:background_gen}, and the estimated covariance assumed in the objective \eqref{eq:rvda_cost} by $\hat B_t$. The misspecified analysis thus replaces $B_t$ with $\hat B_t\neq B_t$ in \eqref{eq:rvda_cost}. 

The misspecification of the background-error covariance may lie in the amplitude of the covariance matrix or in its correlation structure. The amplitude is comparatively straightforward to estimate from innovation statistics \cite{desroziers2005diagnosis, chapnik2006diagnosis}, whereas in practice the correlation structure is hard to specify and is commonly approximated by a simpler one, such as a diagonal covariance \cite{bannister2008review2, stewart2013data}. We therefore hold the amplitude fixed for both $B_t$ and $\hat B_t$ and misspecify only the correlation structure, isolating the impact of an incorrect structural assumption. The same type of misspecification has been studied for the observation-error covariance, where a correlated error is assimilated under a diagonal assumption \cite{stewart2013data}. In our investigation, we apply it instead to the background-error covariance.

We define both covariances in the separable form
\begin{equation}\label{eq:B_decomp}
B_t = \sigma_b^2C, \quad \hat{B}_t = \sigma_b^2\hat{C},
\end{equation}
where $\sigma_b^2$ is the background-error variance and $C$, $\hat C$ are correlation matrices that encode the spatial correlation structure. This separation lets us treat the amplitude and the correlation structure independently. Since the amplitude is readily estimated, we hold $\sigma_b^2$ fixed for both covariances and vary only the correlation structure, so that the misspecification is confined to the discrepancy between $C$ and $\hat C$, with $\hat C = C$ corresponding to the self-consistent case and  $\hat C\neq C$ to the misspecified one.

Having established the roles of the true and assumed covariances and their separation into amplitude and correlation structure, we now specify the correlation matrices themselves. We consider two structures. The first is diagonal, with the correlation matrix equal to the identity, so that the background error is spatially uncorrelated. The second is the correlated structure given by the second-order autoregressive (SOAR) kernel
\begin{equation}\label{eq:soar}
C(r) = \left(1 + \frac{r}{L}\right) e^{-r/L},
\end{equation}
where $r$ is the distance between two grid points and $L$ is the correlation length, controlling the spatial extent over which the background error remains correlated. In our experiments we set $\sigma_b^2=1$ and $L=3\Delta s$. A small ridge is added to the SOAR covariance, giving $B_t=\sigma_b^2C+\epsilon I$ with $\epsilon=10^{-10}$, to ensure numerical positive definiteness.
The SOAR kernel \cite{balgovind1983stochastic, daley1991atmospheric} is a common choice for error correlations in variational data assimilation \cite{ingleby2001statistical, haben2011conditioning, stewart2013data, ebtehaj2014variational}.

Note that the misspecification is not restricted to these two structures. We choose them because their pronounced difference, an uncorrelated structure against a smoothly correlated one, makes the impact of the mismatch most visible.

\section{Numerical Results}
\label{sec:numerics}

Our numerical experiments are designed to investigate the effect of misspecifying the background-error covariance between its generating and assumed roles. We proceed through three settings of increasing difficulty. The first (and simplest) case  considers a self-consistent diagonal covariance  in which the background error is assumed spatially uncorrelated. The second is a self-consistent SOAR covariance defined in \eqref{eq:soar}, representing the best case attainable when the true correlation structure is more complex than a diagonal one yet is correctly specified in the objective. The third is the misspecified case, in which the background error is generated from a SOAR-correlated covariance but assumed diagonal in the objective. We evaluate all three settings on four conservation-law test problems. For the misspecified case, we further examine the effect across multiple assimilation cycles, and show that the regularized variational method with the residual prior transform is more robust than both the Gaussian baseline and the standard total variation prior.

\subsection{Experiment setup}
\label{sec:setup}
All four problems are one-dimensional and develop shocks or other sharp features, either forming from smooth initial data through nonlinear steepening or present from discontinuous initial data, which makes them natural test cases for the structure-preserving reconstruction.

\medskip
\noindent\textbf{Inviscid Burgers' equation on $x \in [-\pi, \pi]$.}
    \begin{equation}\label{eq:burgers}
        u_t + \left(\tfrac{1}{2}u^2\right)_x = 0,
    \end{equation}
    with  $u(x, 0) = 1 + \tfrac{1}{2}\sin x$ and periodic boundary conditions. Since $u(x,0)$ is smooth a shock forms after the wave-breaking time $t_b = 2$. We assimilate past shock formation at $t_f = \pi$ on a grid of $n = 64$ cells.

\medskip
    \noindent\textbf{Stoker problem on $x \in [-1, 1]$}.  
    The governing shallow water equations are given by
    \begin{equation}\label{eq:swe}
        \begin{aligned}
            h_t + (hu)_x &= 0, \\
            (hu)_t + \left(hu^2 + \tfrac{1}{2}g h^2\right)_x &= 0,
        \end{aligned}
    \end{equation}
    where $h$ is the water depth, $u$ the velocity, and  $g = 9.81$ denotes the acceleration due to gravity. We consider the dam-break initial condition \cite{stoker1958water} specified as
    \begin{equation}\label{eq:stoker_ic}
        (h, u)(x, 0) =
        \begin{cases}
            (1, \, 0), & x < 0, \\
            (0.2, \, 0), & x \ge 0,
        \end{cases}
    \end{equation}
    which we assimilate at $t_f = 0.15$ on a grid of $n = 128$ cells.

\medskip
    \noindent\textbf{Sod shock tube \cite{sod1978survey} on $x \in [0, 1]$}.  The governing Euler equations take the form
    \begin{equation}\label{eq:euler}
        \begin{aligned}
            \rho_t + (\rho u)_x &= 0, \\
            (\rho u)_t + (\rho u^2 + p)_x &= 0, \\
            E_t + \bigl((E + p) u\bigr)_x &= 0,
        \end{aligned}
    \end{equation}
    with total energy $E = \frac{p}{\gamma - 1} + \tfrac{1}{2}\rho u^2$ and ratio of specific heats $\gamma = 1.4$. Here $\rho$ is the density, $u$ the velocity, and $p$ the pressure, and the state is assimilated in the primitive variables $(p, u, T)$ with temperature $T = p / \rho$. The Sod initial condition \cite{shu1988efficient} is given by
    \begin{equation}\label{eq:sod_ic}
        (\rho, u, p)(x, 0) =
        \begin{cases}
            (1, \, 0, \, 1), & x \le 0.5, \\
            (0.125, \, 0, \, 0.1), & x > 0.5,
        \end{cases}
    \end{equation}
    which we assimilate at $t_f = 0.2$ on a grid of $n = 128$ cells.

\medskip
    \noindent\textbf{Shu-Osher problem on $x \in [-5, 5]$}. Euler equations \eqref{eq:euler} are used to describe a Mach-3 shock interacting with a sinusoidal density field. The initial condition is prescribed by \cite{shu1988efficient}
    \begin{equation}\label{eq:shu_osher_ic}
        (\rho, u, p)(x, 0) =
        \begin{cases}
            (3.857143, \, 2.629369, \, 10.33333), & x \le -4, \\
            (1 + 0.2\sin(5x), \, 0, \, 1), & x > -4,
        \end{cases}
    \end{equation}
    which we assimilate at $t_f = 1.8$ on a grid of $n = 512$ cells.

The grid resolution $n$ is chosen according to the complexity of each solution profile. In particular, the Shu-Osher problem uses $n=512$ to resolve its fine-scale post-shock structure, which is essential for a meaningful comparison.
For the Stoker and Sod problems, the true profile is available in closed form as the exact Riemann solution \cite{stoker1958water, sod1978survey}, which we cell-average onto the assimilation grid. No closed-form solution is available for either the Burgers or Shu-Osher problem at the final time.  In these cases we use highly resolved numerical approximations, respectively with $n = 4096$ and $n=10240$, and project onto the assimilation grid.  We employ the fifth-order WENO \cite{jiang1996efficient}  coupled with Lax-Friedrichs flux splitting for the Burgers and Stoker problems, and an HLLC finite volume scheme \cite{toro1994hllc} with MUSCL reconstruction \cite{van1979towards} for both Euler problems.
In all cases the time integration uses a third-order strong-stability-preserving Runge-Kutta method \cite{gottlieb1998total}.

The choice of numerical schemes is intentional and reflects a trade-off inherent to hyperbolic solvers. Specifically, some numerical dissipation is required for stability, but excessive dissipation prevents the solver from preserving structure once shocks form. This trade-off is compounded in a data assimilation setting, since 3D-Var feeds external information into the model at each analysis step, and any non-physical oscillatory behavior propagating through the assimilation cycle may be amplified, potentially steering successive analyses away from the true state rather than merely degrading a single forecast.

The adaptive weighting of (5th order) WENO is known to be effective at preserving structure for solving hyperbolic PDEs that admit shock discontinuities, and we use it here for forecasting Burgers' equation and Stoker's problem. While some non-physical, method-induced oscillations are present, they do not cause any long-term instabilities in the 3DVar solution.  The presence of strong shocks in the Euler equations (Sod and Shu–Osher) makes oscillations of this same kind a more serious concern, particularly where they interact with a fine-scale entropy wave, as in the Shu–Osher case.   We therefore adopt HLLC-MUSCL for solving the Euler equations. Although more dissipative than WENO, its reliance on a TVD limiter means that the solution remains robust in these regimes \cite{toro2009riemann, sweby1984tvd}. The scheme resolves the contact and shock waves without spurious oscillations because HLLC's approximate Riemann solver explicitly models the wave structure at the flux level \cite{toro1994hllc}, making a componentwise reconstruction sufficient for this purpose.

\subsubsection*{Observation setup} The observation-error covariance in \eqref{eq:observation_model} takes the simplest spatially uncorrelated structure, so that $\Sigma_t$ is diagonal. The noise level is set by the signal-to-noise ratio (SNR), specified in decibels. For each observed variable with true values $\bm f$ at the observed cells, we define the SNR as
\begin{equation}\label{eq:snr}
\text{SNR} = 10\log_{10}\frac{\sigma_f^2}{\sigma_o^2}, \quad \sigma_f^2 = \frac{1}{m}\sum_{i=1}^{m} (f_i - \bar{f})^2,
\end{equation}
where $\sigma_f$ is the root-mean-square amplitude of the signal about its mean $\bar f$, and $\sigma_o$ the observation-error standard deviation.

We prescribe SNR $=10$ dB for all tests except when specifically performing SNR analysis  (see \Cref{fig:sweep_snr}). Moreover, in all cases the observation matrix $H$ in \eqref{eq:observation_model} subsamples the state with stride $3$, so that observations are available at every third cell.  The  exception is when we perform  observation sparsity analysis (see \Cref{fig:sweep_stride}). 

\subsubsection*{Performance metrics} Each experiment is repeated over $100$ Monte Carlo realizations of the observation noise, from which we report the ensemble mean of each method. 
To ensure a meaningful comparison, we display a method's mean reconstruction only when at least 10 realizations yield a physically valid state.
In each experiment we report the pointwise error
\begin{equation}\label{eq:pw_error}
E_j^{abs} = \abs{\hat{f}_j - f_j},
\end{equation}
where $\hat f_j$ is the $j$th component of the mean analysis over the valid realizations and $f_j$ the corresponding value of the true PDE profile. 

In each figure we denote as $G$  the reconstruction obtained from \eqref{eq:3dvar_cost} {\em without} regularization, which serves as the unregularized baseline. The solutions labeled $TV$ and $R$ are obtained from \eqref{eq:rvda_cost} with the standard (first order) total variation (TV) prior transform and the residual prior transform, respectively. The TV prior transform serves as the standard structure-preserving benchmark against which we measure the residual prior transform. Under this choice the total variation prior imposes a fixed piecewise constant assumption, whereas the residual prior transform does not commit to any fixed order of smoothness. The comparison therefore tests how each prior performs when the structure of the solution is unknown in advance, where we expect the residual prior transform to be more robust.  
For ease of presentation we display only the main variable of each problem: $u$ for the Burgers problem, $h$ for the Stoker problem, and the pressure $p$ for the Sod and Shu-Osher problems, taking a native state variable in each case.

\subsection{Self-consistent background covariance}
\label{sec:selfconsist}

\begin{figure}[!ht] 
    \centering
    \begin{subfigure}[b]{0.24\textwidth}
        \includegraphics[width=\textwidth]{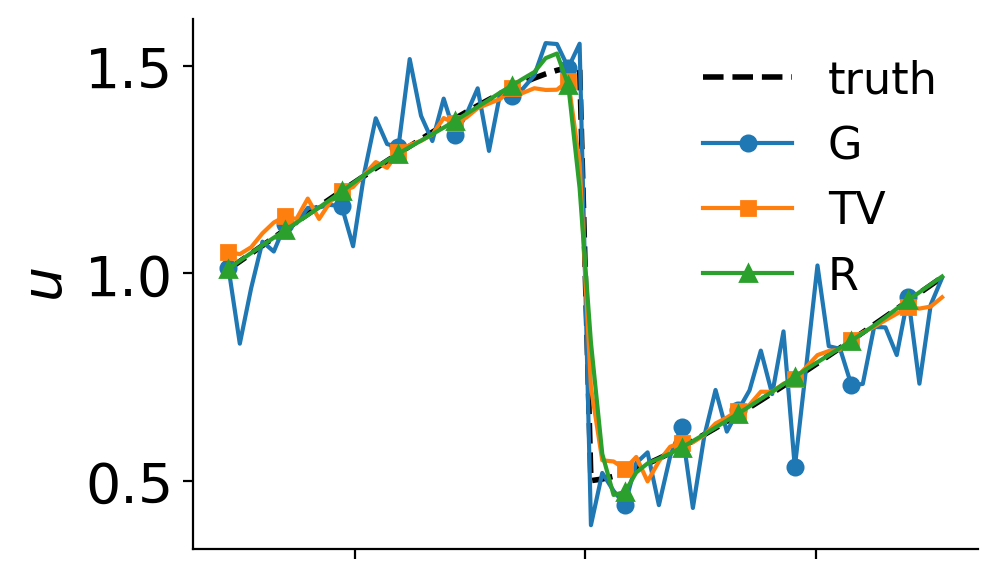}
    \end{subfigure}
    \hfill
    \begin{subfigure}[b]{0.24\textwidth}
        \includegraphics[width=\textwidth]{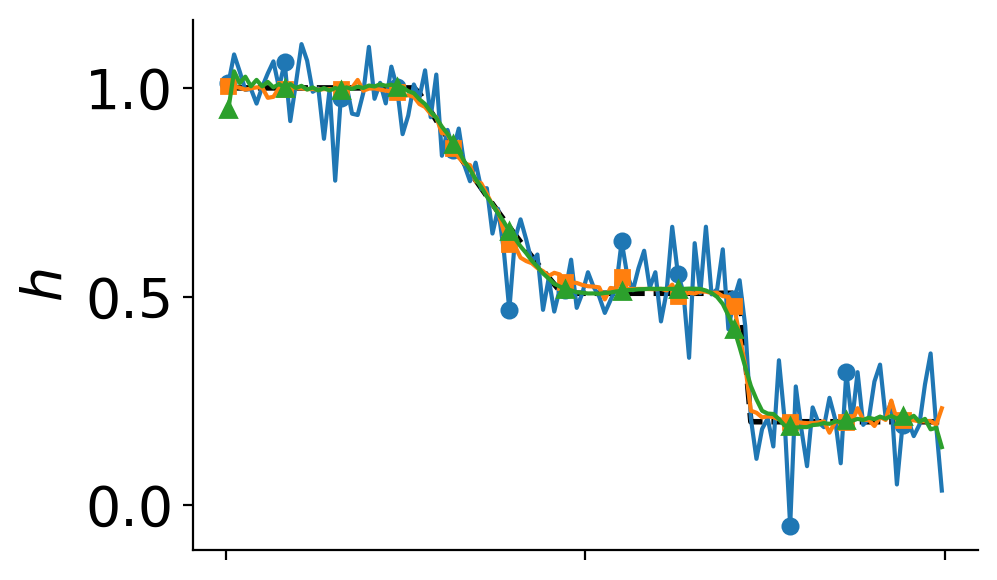}
    \end{subfigure}
    \hfill
    \begin{subfigure}[b]{0.24\textwidth}
        \includegraphics[width=\textwidth]{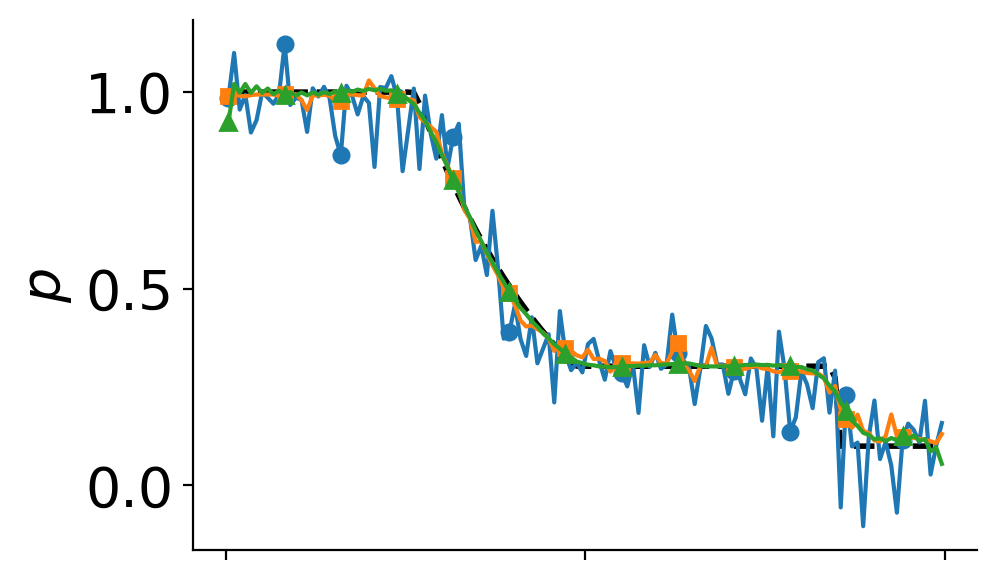}
    \end{subfigure}
    \hfill
    \begin{subfigure}[b]{0.24\textwidth}
        \includegraphics[width=\textwidth]{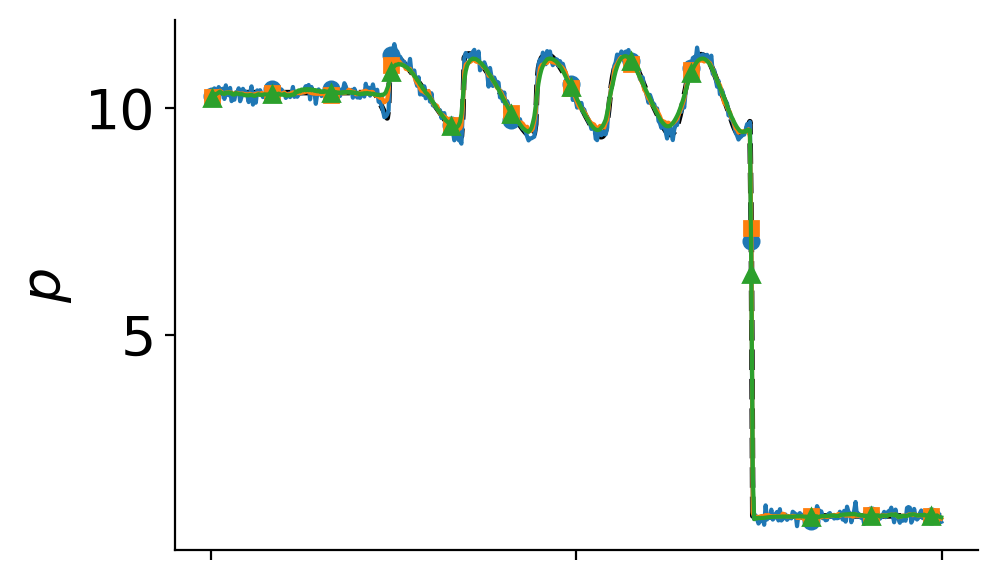}
    \end{subfigure}

    \vspace{0.5em}

    \begin{subfigure}[b]{0.24\textwidth}
        \includegraphics[width=\textwidth]{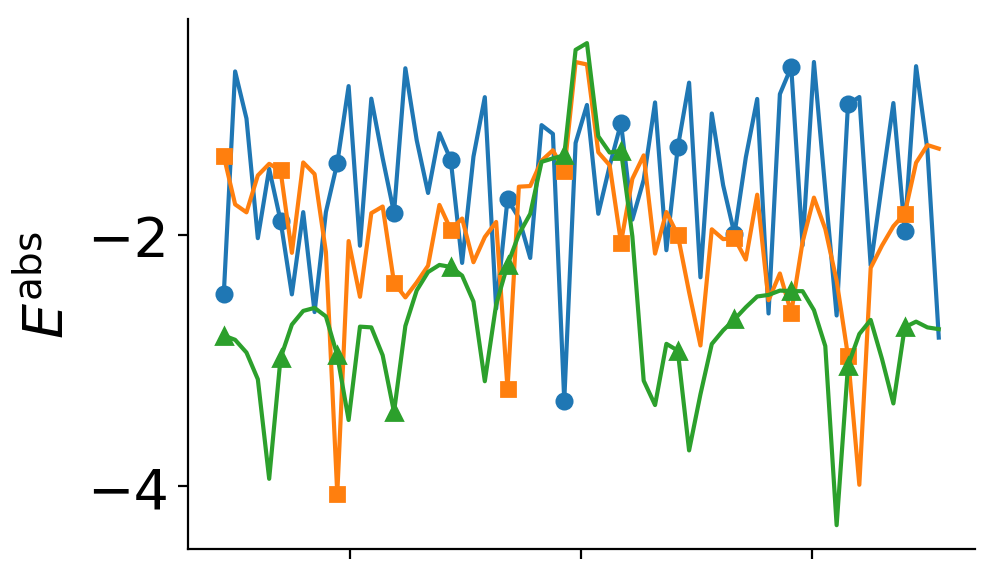}
    \end{subfigure}
    \hfill
    \begin{subfigure}[b]{0.24\textwidth}
        \includegraphics[width=\textwidth]{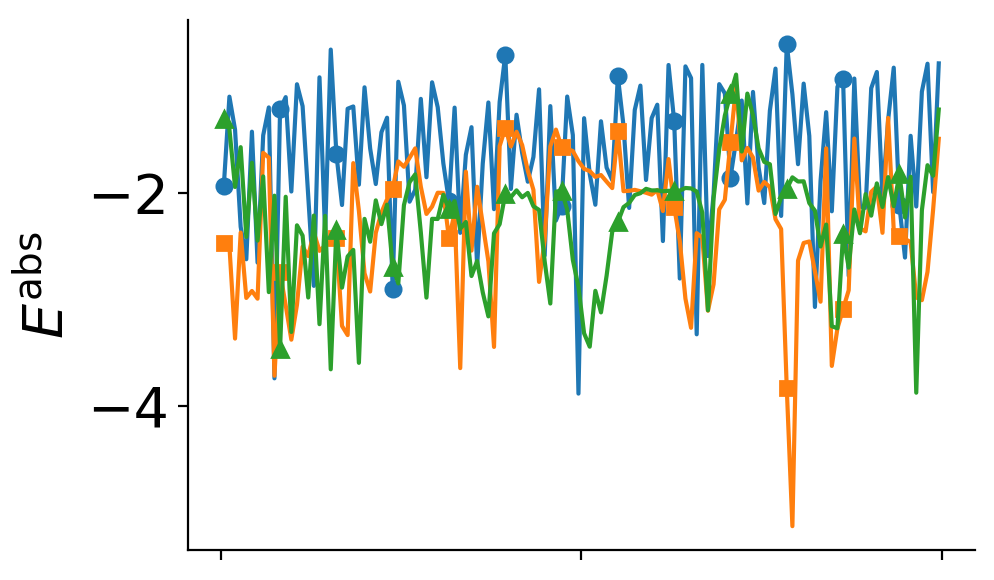}
    \end{subfigure}
    \hfill
    \begin{subfigure}[b]{0.24\textwidth}
        \includegraphics[width=\textwidth]{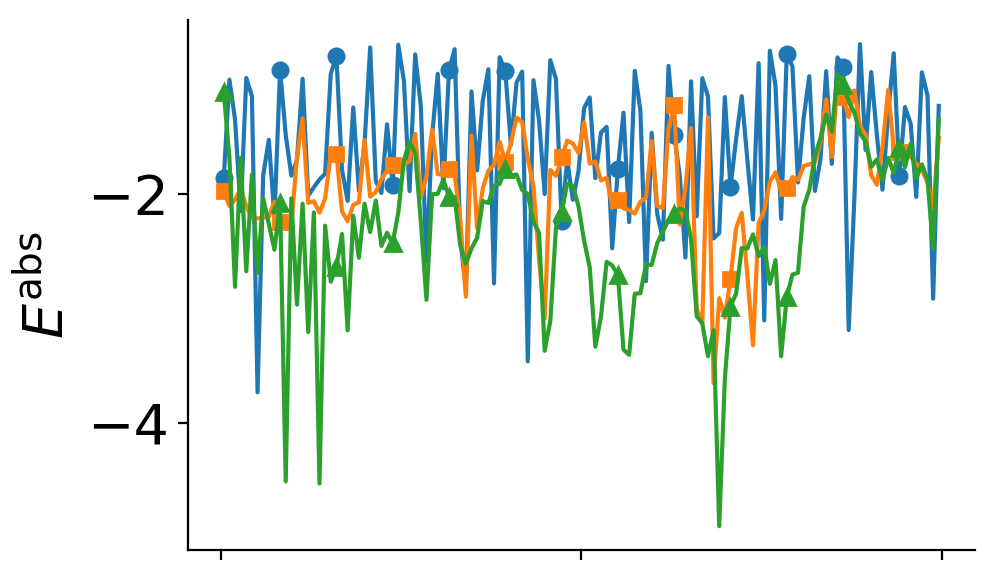}
    \end{subfigure}
    \hfill
    \begin{subfigure}[b]{0.24\textwidth}
        \includegraphics[width=\textwidth]{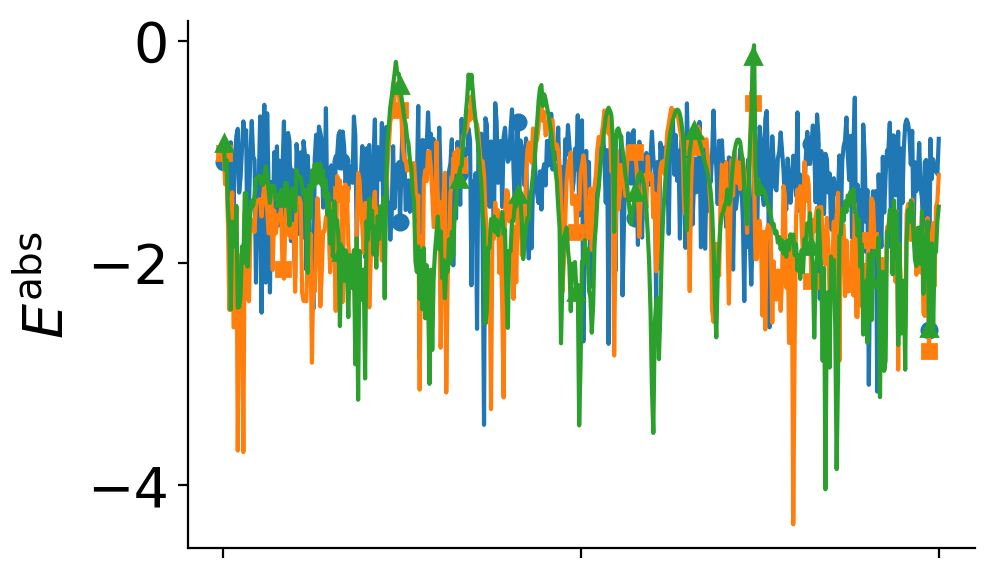}
    \end{subfigure}

    \caption{The background-error covariance is diagonal-structured. (top) Truth profile and mean reconstructions of the representative variable obtained by Gaussian, TV prior transform, and residual prior transform for (left) Burgers, (middle left) Stoker, (middle right) Sod, and (right) Shu-Osher. (bottom) The corresponding pointwise error \eqref{eq:pw_error} on a $\log_{10}$ scale.}
    \label{fig:selfconsist_diag_single}
\end{figure}

Our investigation starts with a diagonal background-error covariance that is self-consistent between the generation in \eqref{eq:background_gen} and the objective in \eqref{eq:rvda_cost}, in the single-cycle setting with $T=1$ in \eqref{eq:forecast_model}. This diagonal case is the simplest of the self-consistent settings and verifies the methods under ideal conditions. We then consider a SOAR-structured background-error covariance, still self-consistent, which represents the best case each method can achieve when the true correlation structure is more complex than a diagonal one yet is correctly specified in the objective.

\begin{figure}[!ht]
    \centering
    \begin{subfigure}[b]{0.24\textwidth}
        \includegraphics[width=\textwidth]{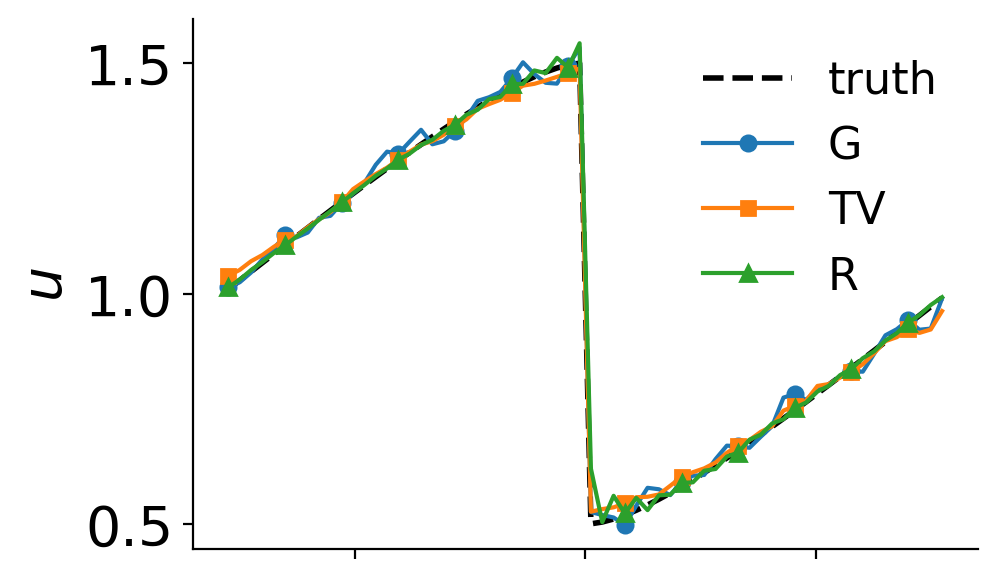}
    \end{subfigure}
    \hfill
    \begin{subfigure}[b]{0.24\textwidth}
        \includegraphics[width=\textwidth]{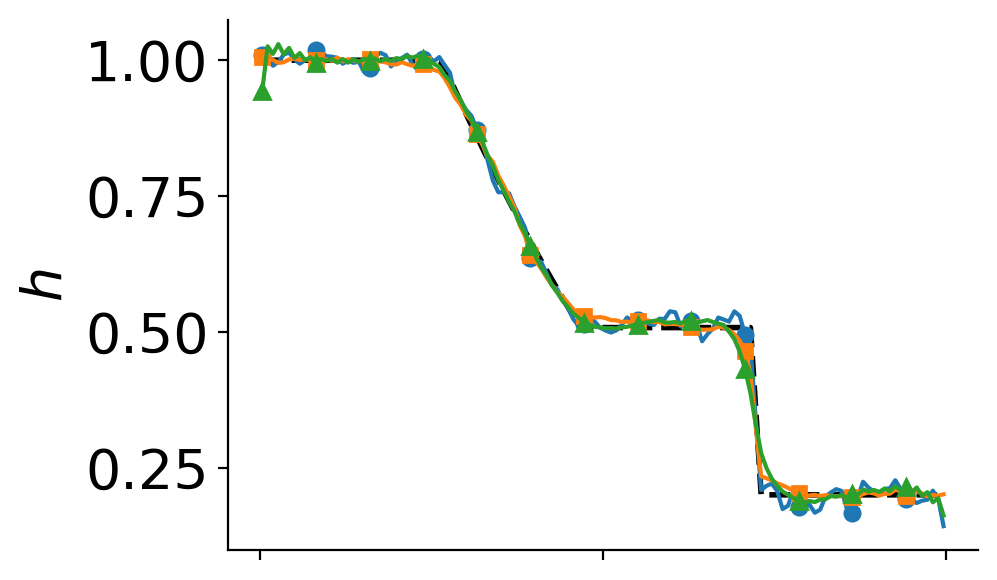}
    \end{subfigure}
    \hfill
    \begin{subfigure}[b]{0.24\textwidth}
        \includegraphics[width=\textwidth]{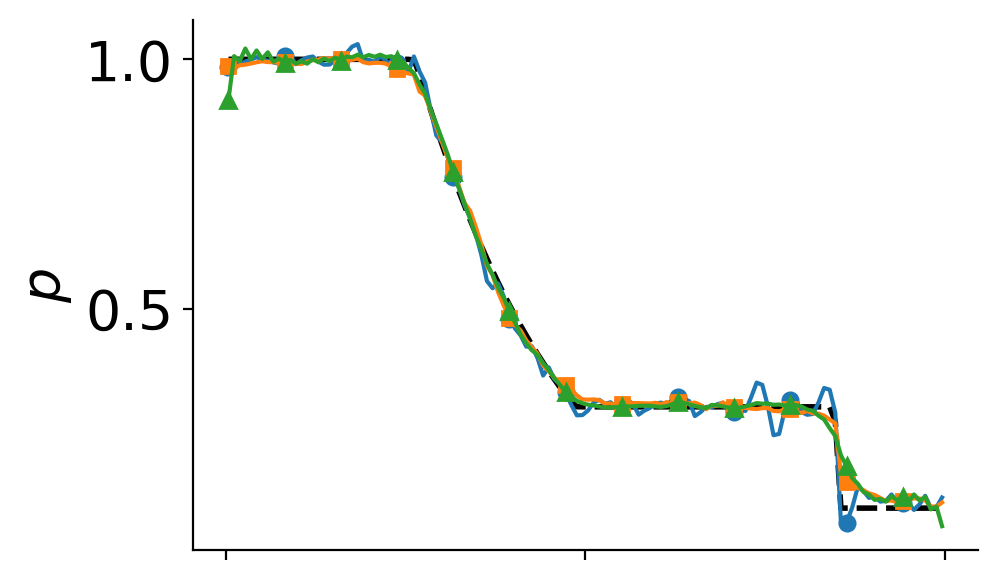}
    \end{subfigure}
    \hfill
    \begin{subfigure}[b]{0.24\textwidth}
        \includegraphics[width=\textwidth]{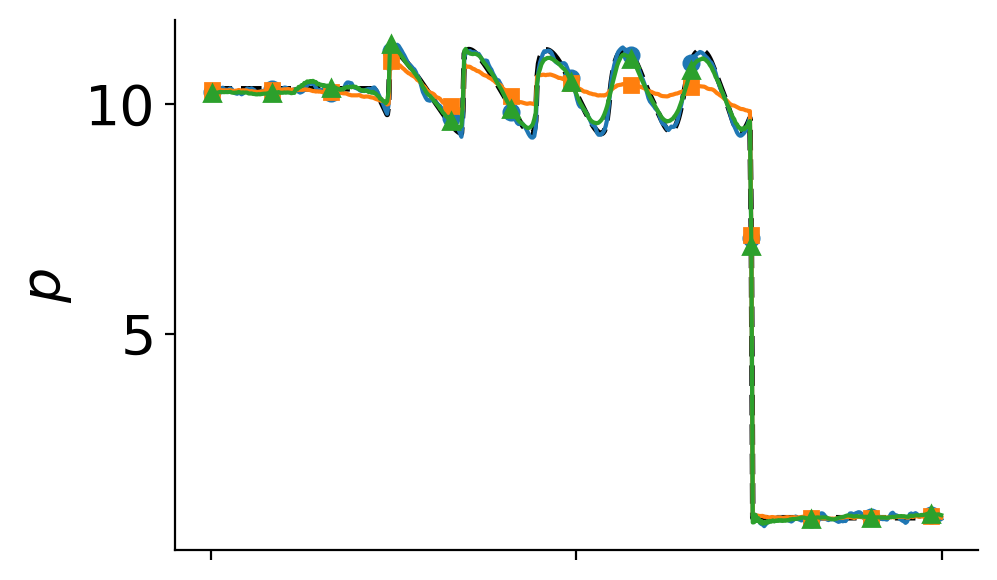}
    \end{subfigure}

    \vspace{0.5em}

    \begin{subfigure}[b]{0.24\textwidth}
        \includegraphics[width=\textwidth]{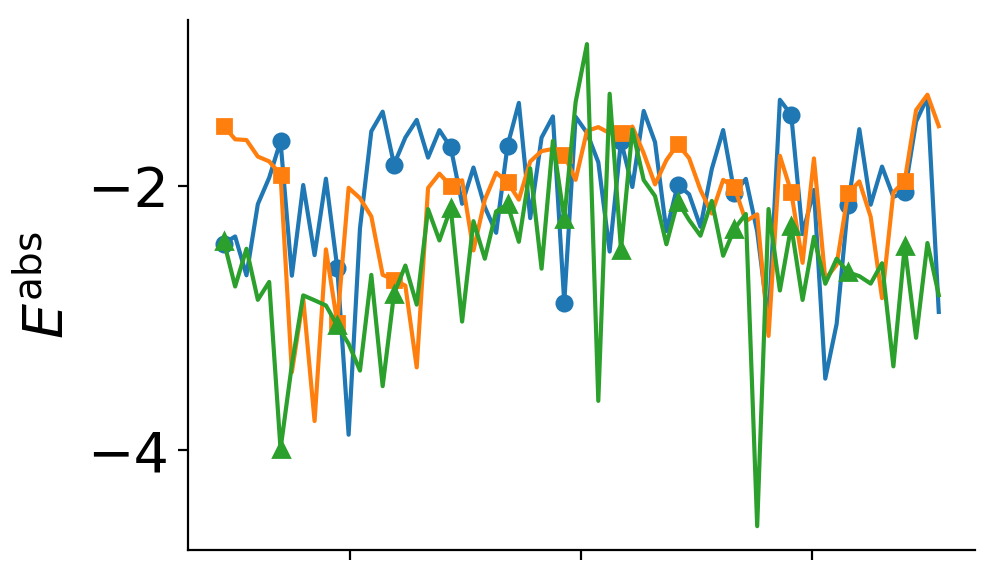}
    \end{subfigure}
    \hfill
    \begin{subfigure}[b]{0.24\textwidth}
        \includegraphics[width=\textwidth]{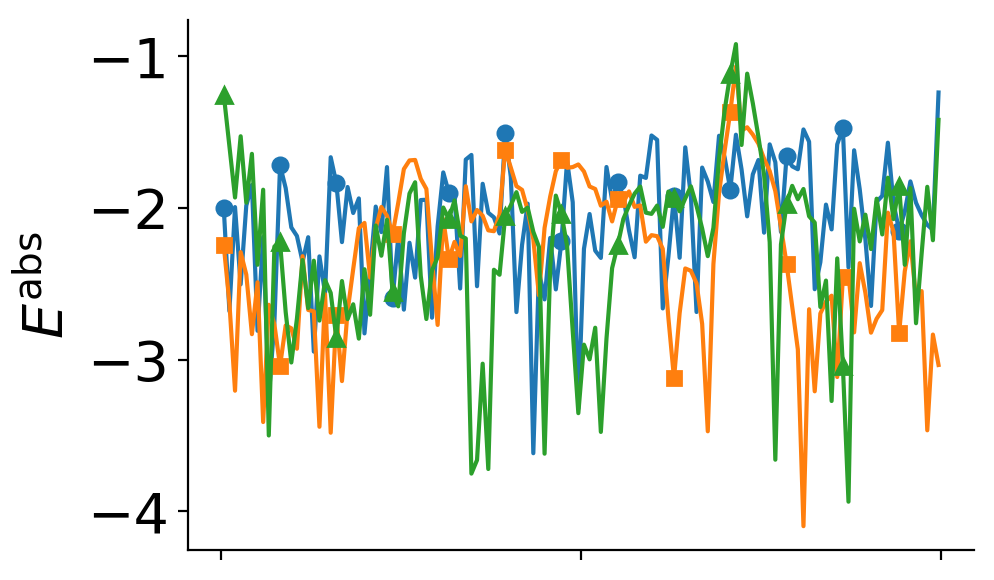}
    \end{subfigure}
    \hfill
    \begin{subfigure}[b]{0.24\textwidth}
        \includegraphics[width=\textwidth]{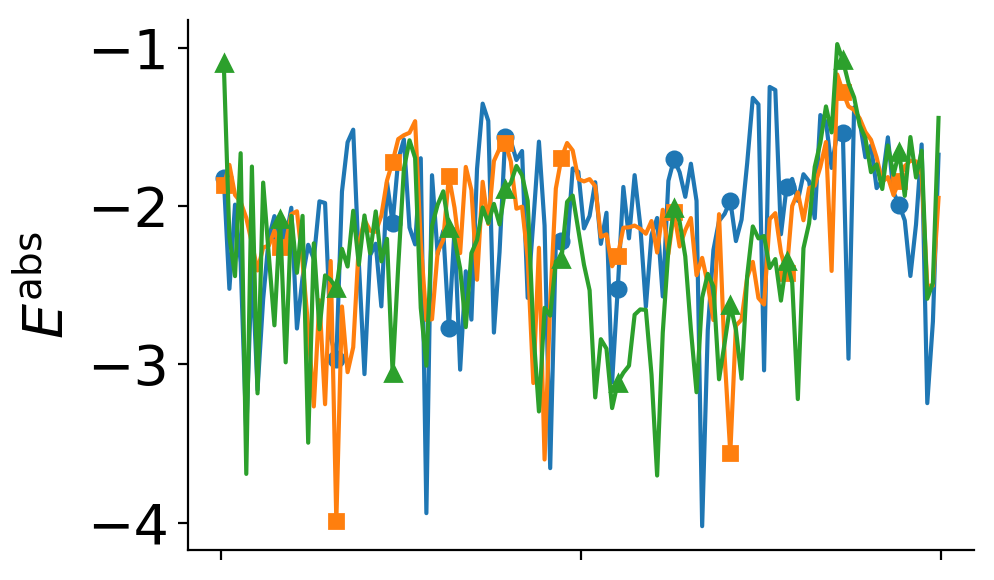}
    \end{subfigure}
    \hfill
    \begin{subfigure}[b]{0.24\textwidth}
        \includegraphics[width=\textwidth]{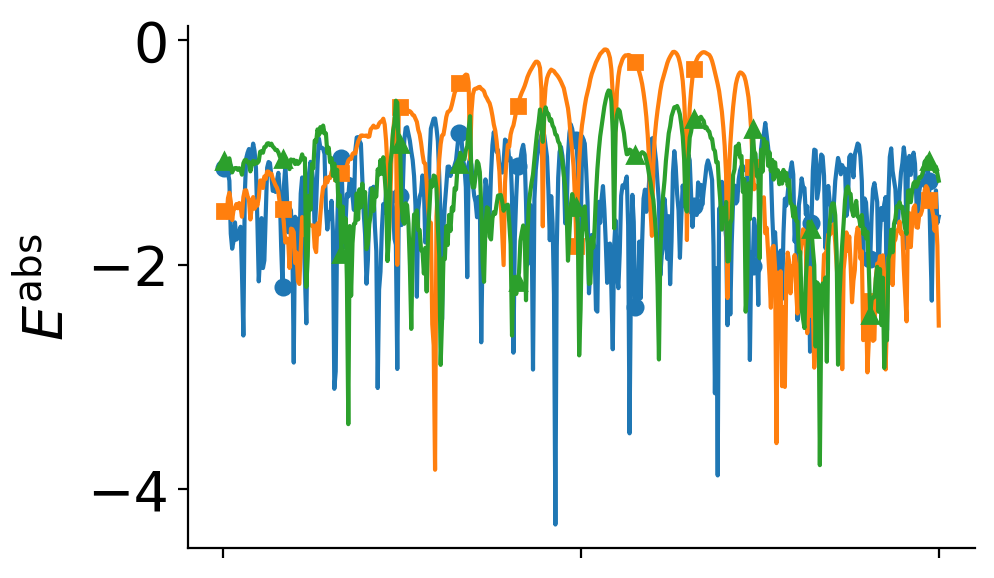}
    \end{subfigure}

    \caption{The background-error covariance is SOAR-correlated and correctly specified in the objective. (top) Truth profile and mean reconstructions of the representative variable obtained by Gaussian, TV prior transform, and residual prior transform for (left) Burgers, (middle left)Stoker, (middle right) Sod, and (right) Shu-Osher. (bottom) The corresponding pointwise error \eqref{eq:pw_error} on a $\log_{10}$ scale.}
    \label{fig:selfconsist_soar_single}
\end{figure}

\Cref{fig:selfconsist_diag_single}(top) shows the mean reconstruction of each method over the $100$ Monte Carlo simulations under correct specification of diagonal background-error covariance. The bottom row shows the corresponding pointwise $\log_{10}$ error. Both prior transforms in \eqref{eq:rvda_cost} suppress the oscillations over the smooth regions far more effectively than the Gaussian reconstruction. The standard total variation operator, however, produces spurious oscillations where the solution varies substantially, whereas the residual prior transform stays close to the true profile. 

\Cref{fig:selfconsist_soar_single}(top) shows the mean reconstruction of each method over the $100$ Monte Carlo simulations, under a correctly specified SOAR-correlated background-error covariance. For the Burgers, Stoker, and Sod problems, whose solutions have a relatively simple structure, all three methods give comparable reconstructions that stay close to the true profile.  The advantage of the residual prior transform emerges for the Shu-Osher problem, where the more complex post-shock oscillations are captured more accurately by the TV prior transform than by the residual prior transform. The TV prior transform over-regularizes the solution, in particular by flattening genuine peaks and troughs.

\subsection{Misspecified background covariance}
\label{sec:mismatch}

\begin{figure}[!ht]
    \centering
    \begin{subfigure}[b]{0.24\textwidth}
        \includegraphics[width=\textwidth]{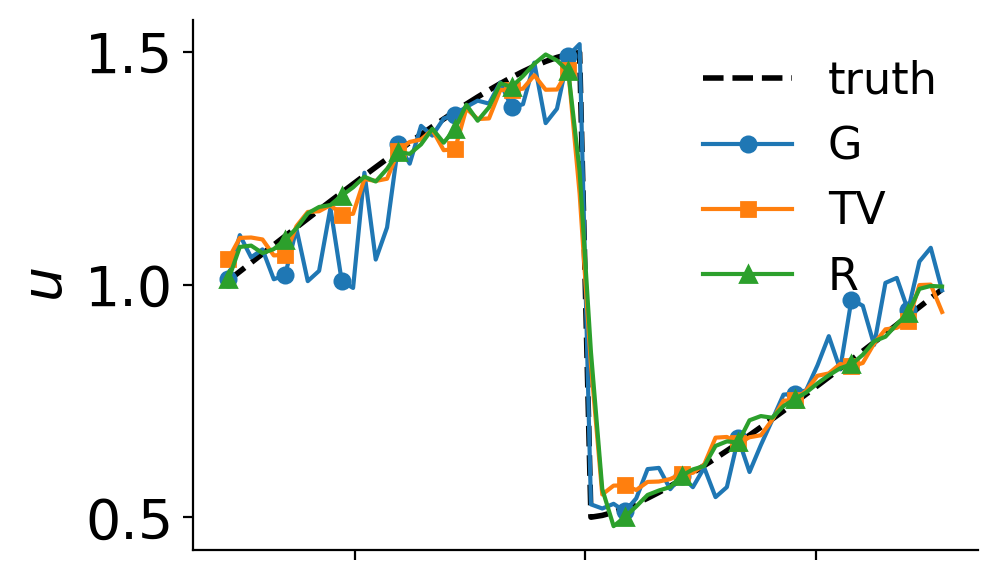}
    \end{subfigure}
    \hfill
    \begin{subfigure}[b]{0.24\textwidth}
        \includegraphics[width=\textwidth]{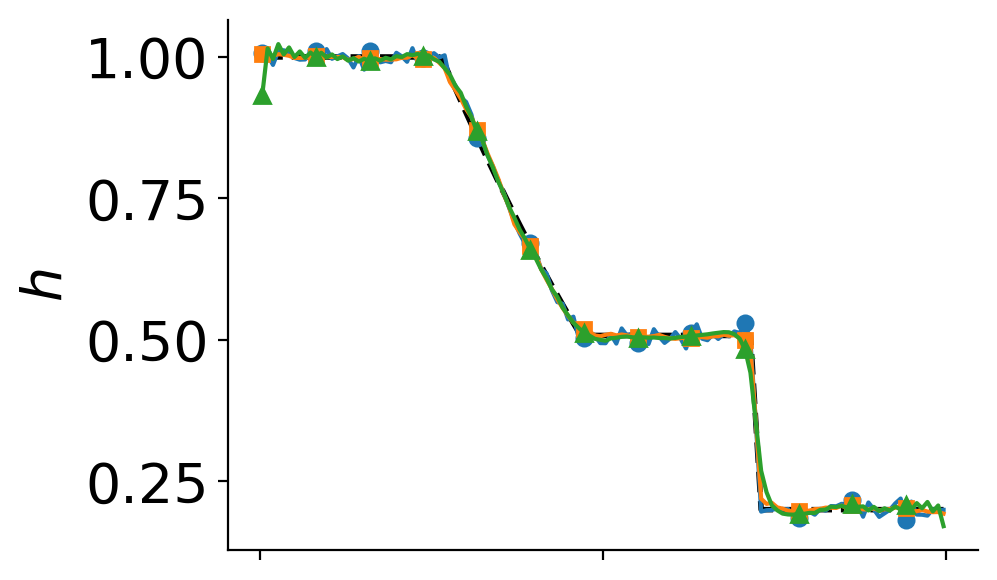}
    \end{subfigure}
    \hfill
    \begin{subfigure}[b]{0.24\textwidth}
        \includegraphics[width=\textwidth]{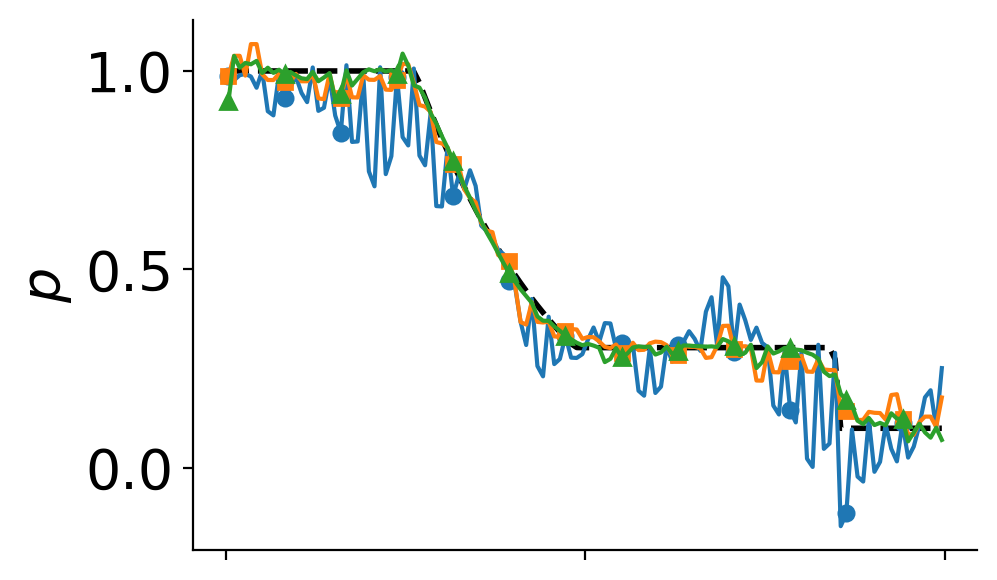}
    \end{subfigure}
    \hfill
    \begin{subfigure}[b]{0.24\textwidth}
        \includegraphics[width=\textwidth]{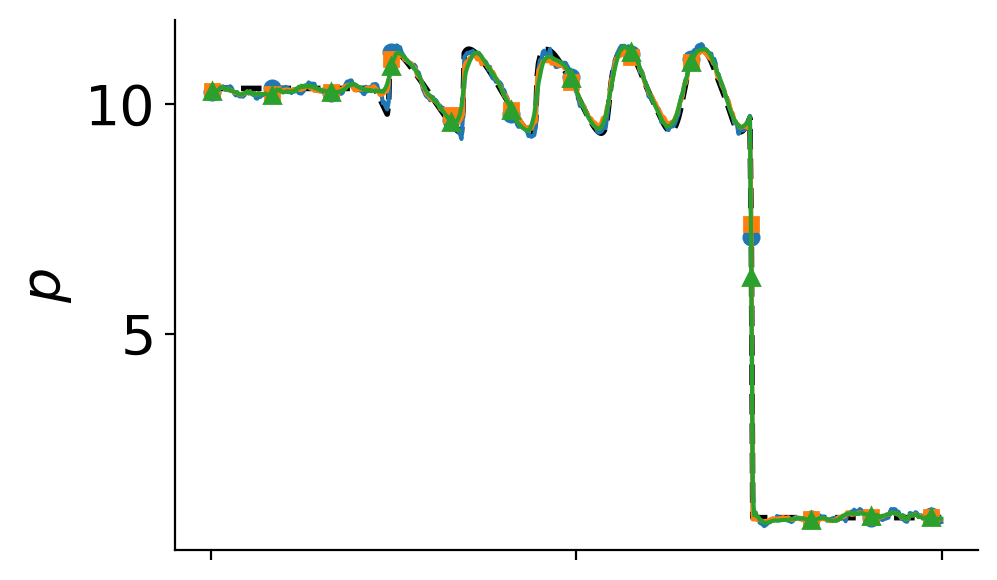}
    \end{subfigure}

    \vspace{0.5em}

    \begin{subfigure}[b]{0.24\textwidth}
        \includegraphics[width=\textwidth]{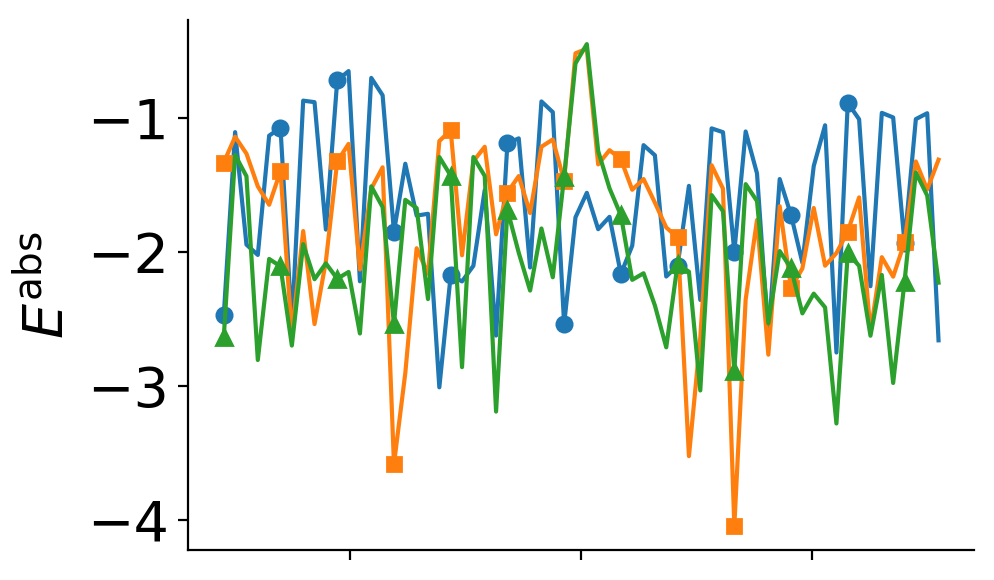}
    \end{subfigure}
    \hfill
    \begin{subfigure}[b]{0.24\textwidth}
        \includegraphics[width=\textwidth]{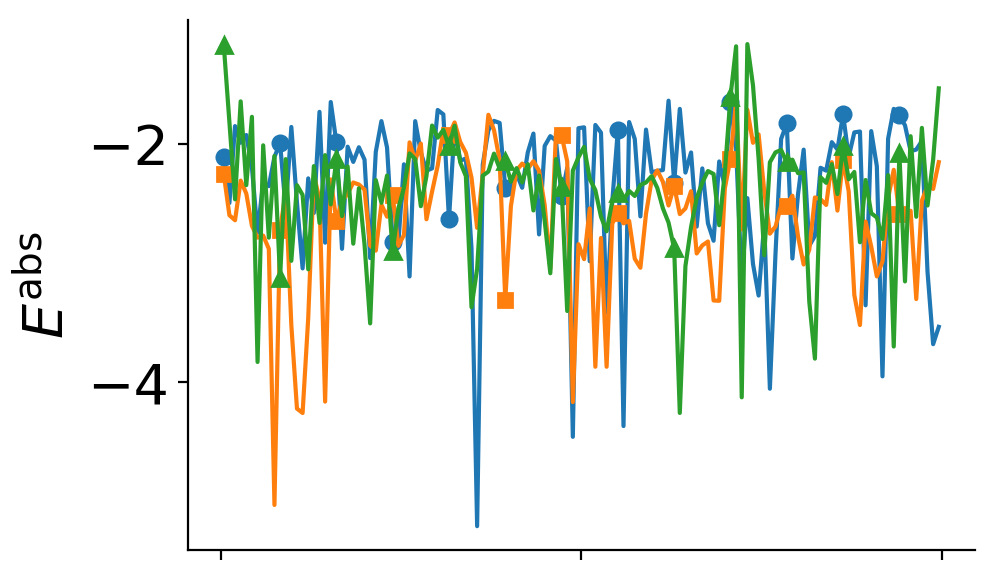}
    \end{subfigure}
    \hfill
    \begin{subfigure}[b]{0.24\textwidth}
        \includegraphics[width=\textwidth]{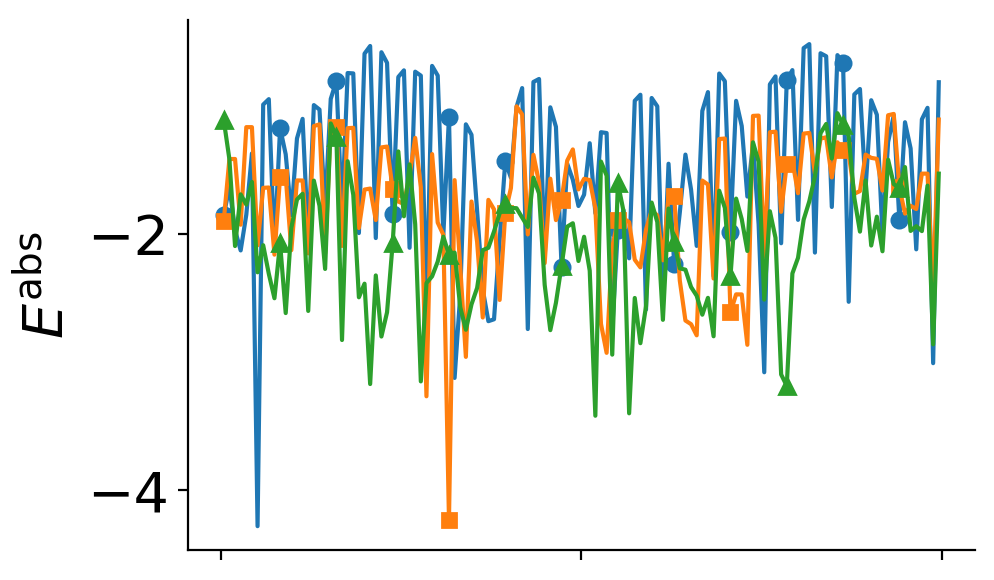}
    \end{subfigure}
    \hfill
    \begin{subfigure}[b]{0.24\textwidth}
        \includegraphics[width=\textwidth]{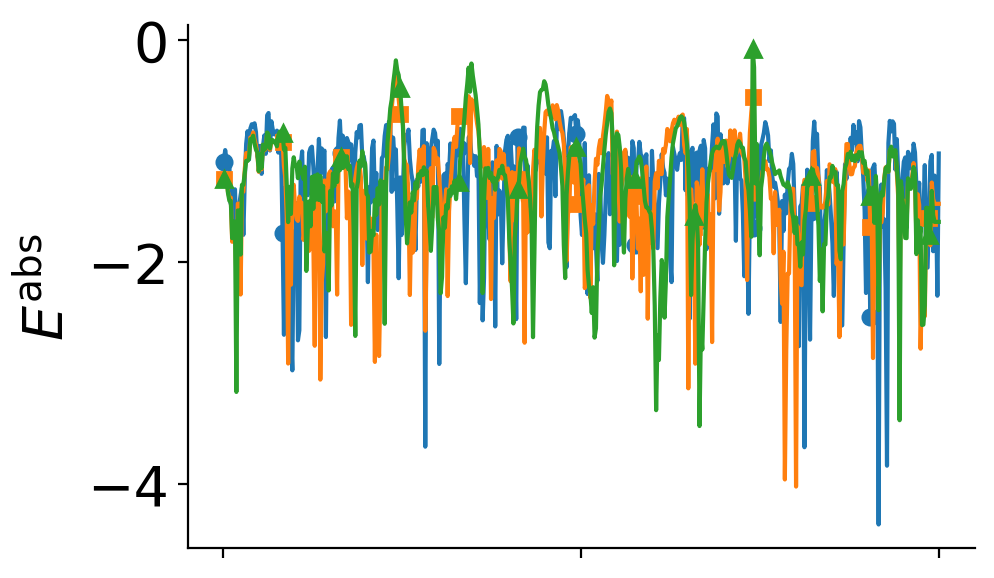}
    \end{subfigure}

    \caption{The background-error covariance is SOAR-correlated in the generation but misspecified as diagonal in the objective, in the single-cycle setting. (top) Truth profile and mean reconstructions of the representative variable obtained by Gaussian, the TV prior transform, and the residual prior transform for (left) Burgers, (middle left) Stoker, (middle right) Sod, and (right) Shu-Osher. (bottom) The corresponding pointwise error \eqref{eq:pw_error} on a $\log_{10}$ scale.}
    \label{fig:mismatch_single}
\end{figure}

We next consider the scenario in which the background-error covariance is SOAR-correlated in the generation of \eqref{eq:background_gen} but misspecified as diagonal in the objective. We evaluate the methods first in the single-cycle setting with $T=1$ in \eqref{eq:forecast_model}, and then in the cycled setting with $T>1$ to assess the stability of the regularized methods over multiple cycles. Since noise enters both the background and the observations at each cycle, causing uncertainty to accumulate through the assimilation process, we expect the results to be worse in the multiple cycle setting.

\begin{figure}[!ht]
    \centering
    \begin{subfigure}[b]{0.32\textwidth}
        \includegraphics[width=\textwidth]{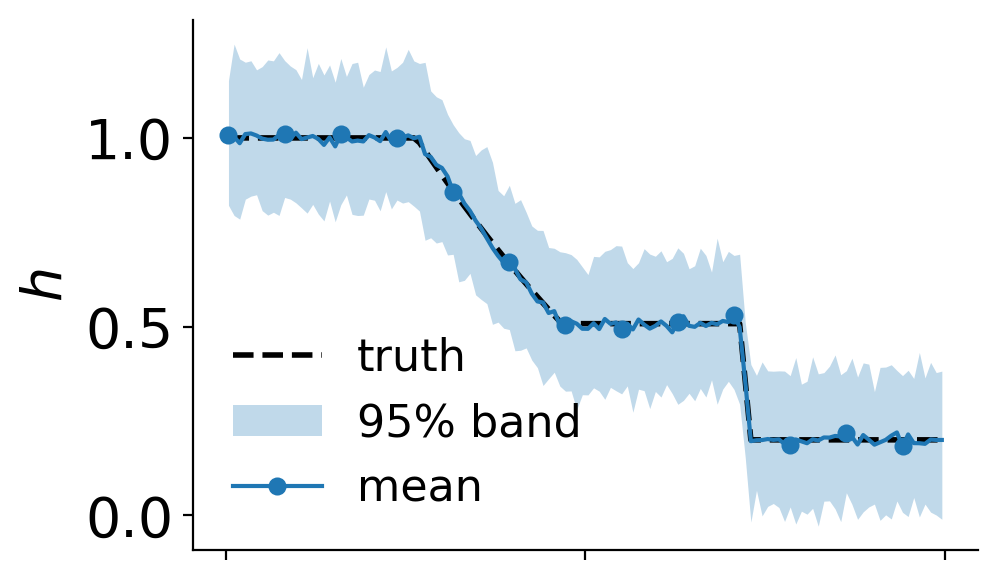}
        \caption{Gaussian}
    \end{subfigure}
    \hfill
    \begin{subfigure}[b]{0.32\textwidth}
        \includegraphics[width=\textwidth]{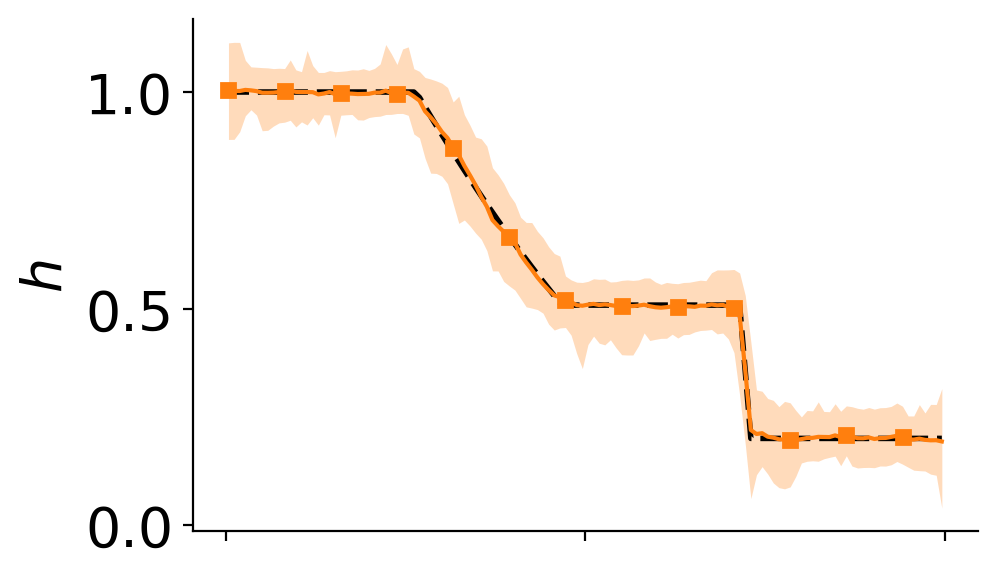}
        \caption{TV prior transform}
    \end{subfigure}
    \hfill
    \begin{subfigure}[b]{0.32\textwidth}
        \includegraphics[width=\textwidth]{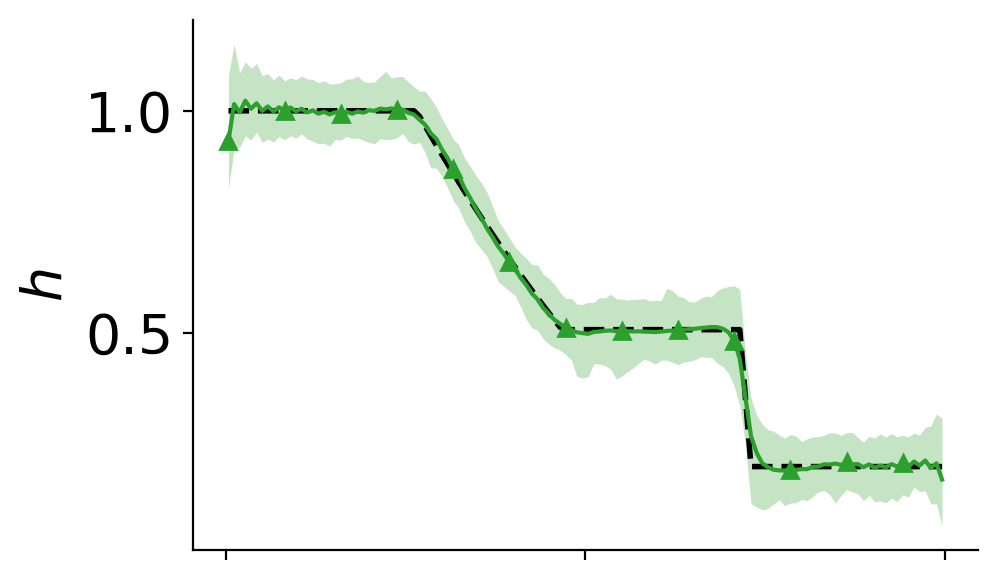}
        \caption{Residual prior transform}
    \end{subfigure}

    \caption{Mean and empirical $95\%$ band of the water depth $h$ for the Stoker problem under the misspecified background-error covariance, single-cycle setting. Each panel shows one method: (left) Gaussian, (middle) total variation, and (right) residual prior transform. The shaded box is the percentile band across the $100$ Monte Carlo realizations.}
    \label{fig:band_mismatch_stoker}
\end{figure}

\Cref{fig:mismatch_single}(top) shows the mean reconstruction of each method over the $100$ Monte Carlo realizations, under a diagonal background-error covariance in the objective while the background is generated with a SOAR-correlated covariance. Because the true correlation structure is more complex than the assumed diagonal one, the standard 3D-Var method \eqref{eq:3dvar_cost} cannot capture the structural behavior, and becomes highly oscillatory over the smooth regions, as in \Cref{fig:selfconsist_diag_single}. Using either the TV or the residual prior transforms in \eqref{eq:rvda_cost} will help to suppress these oscillations and stay close to the true profile overall. The two priors differ across the whole domain, however. For the Burgers and Sod problems, the TV prior introduces spurious oscillations throughout, while the residual prior transform remains stable and smooth.

\begin{figure}[!ht]
    \centering
    \begin{subfigure}[b]{0.24\textwidth}
        \includegraphics[width=\textwidth]{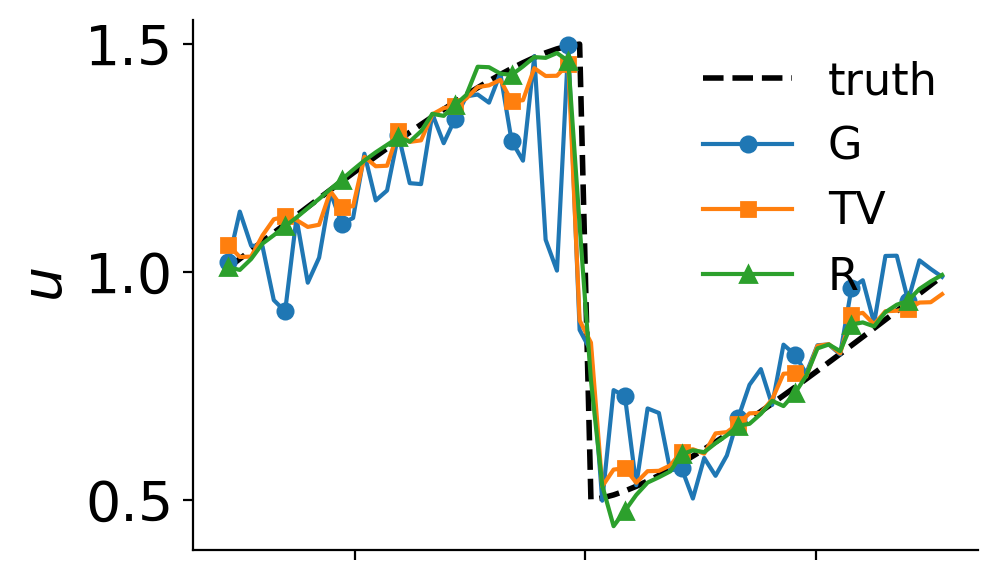}
    \end{subfigure}
    \hfill
    \begin{subfigure}[b]{0.24\textwidth}
        \includegraphics[width=\textwidth]{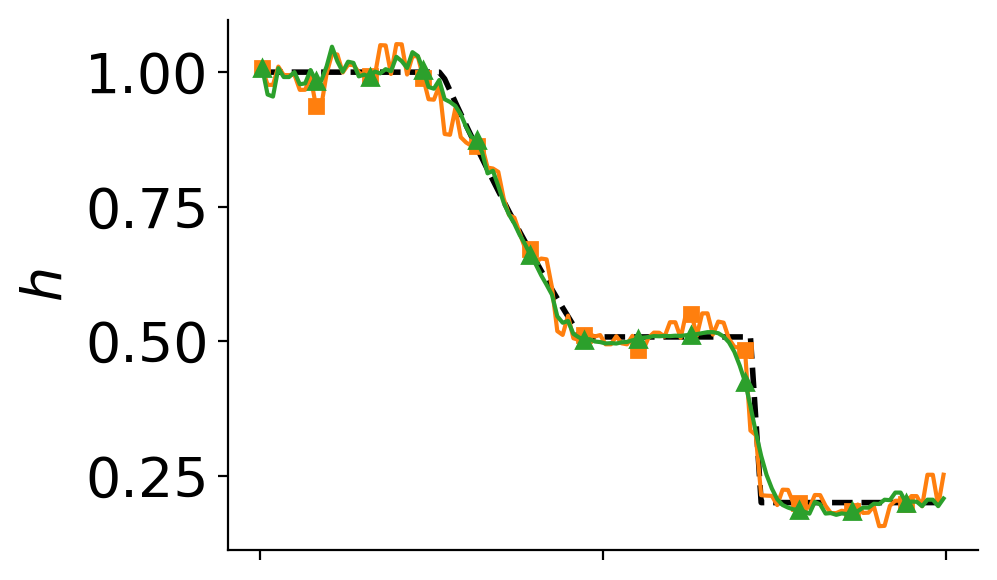}
    \end{subfigure}
    \hfill
    \begin{subfigure}[b]{0.24\textwidth}
        \includegraphics[width=\textwidth]{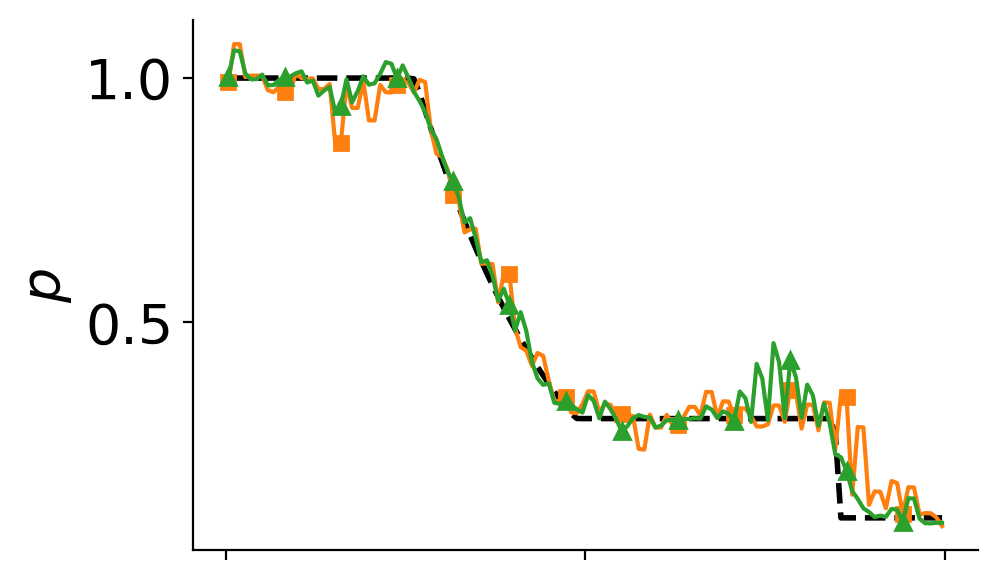}
    \end{subfigure}
    \hfill
    \begin{subfigure}[b]{0.24\textwidth}
        \includegraphics[width=\textwidth]{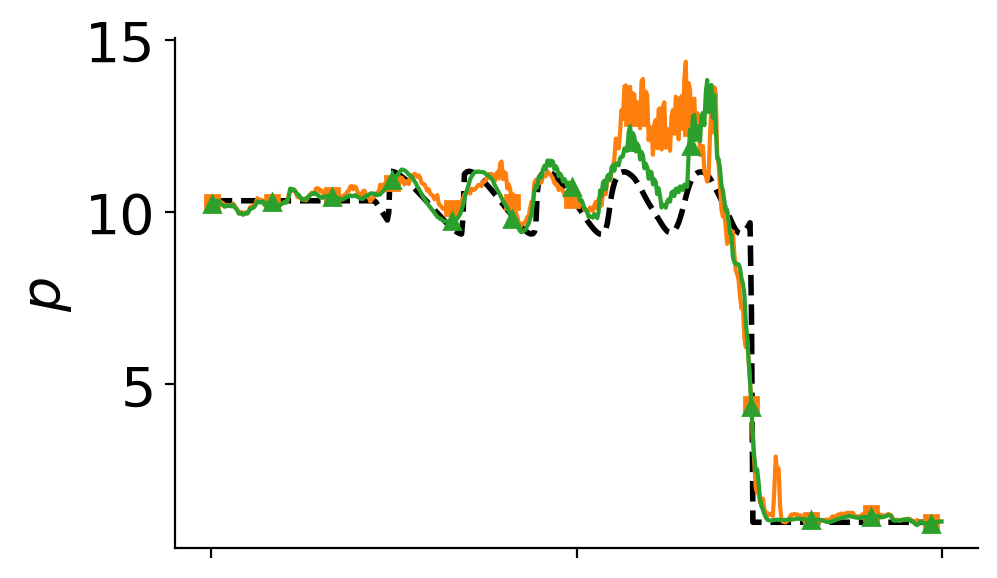}
    \end{subfigure}

    \vspace{0.5em}

    \begin{subfigure}[b]{0.24\textwidth}
        \includegraphics[width=\textwidth]{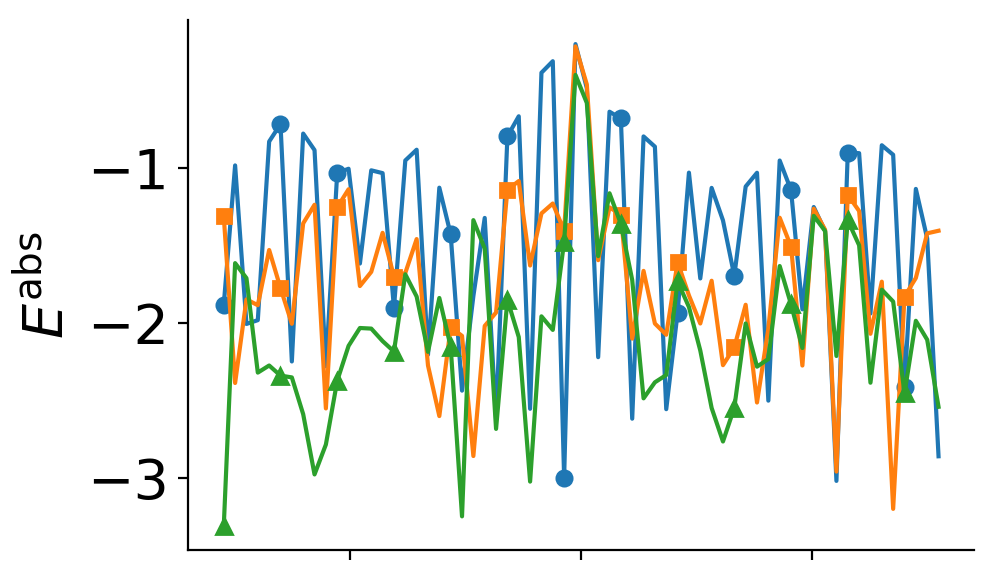}
    \end{subfigure}
    \hfill
    \begin{subfigure}[b]{0.24\textwidth}
        \includegraphics[width=\textwidth]{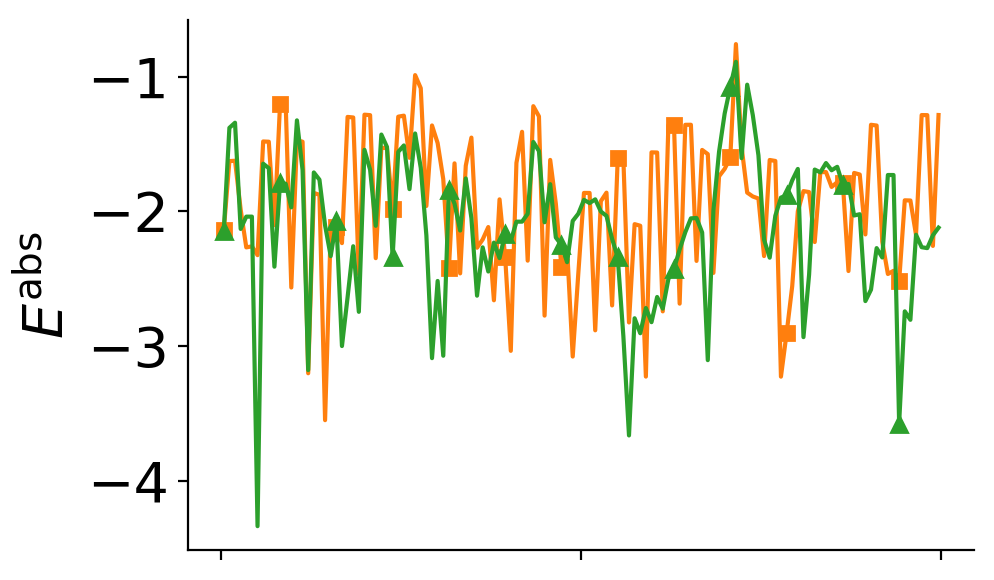}
    \end{subfigure}
    \hfill
    \begin{subfigure}[b]{0.24\textwidth}
        \includegraphics[width=\textwidth]{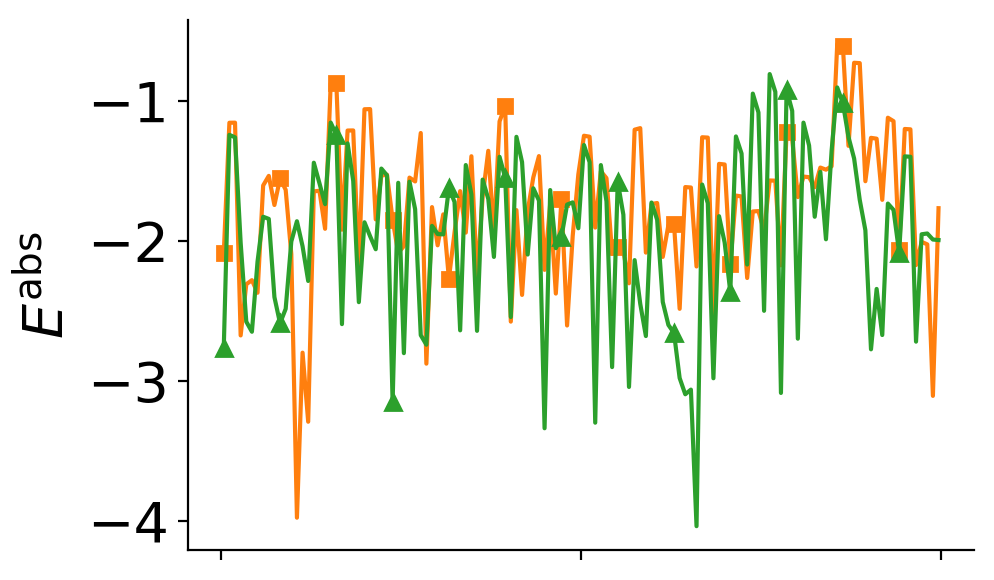}
    \end{subfigure}
    \hfill
    \begin{subfigure}[b]{0.24\textwidth}
        \includegraphics[width=\textwidth]{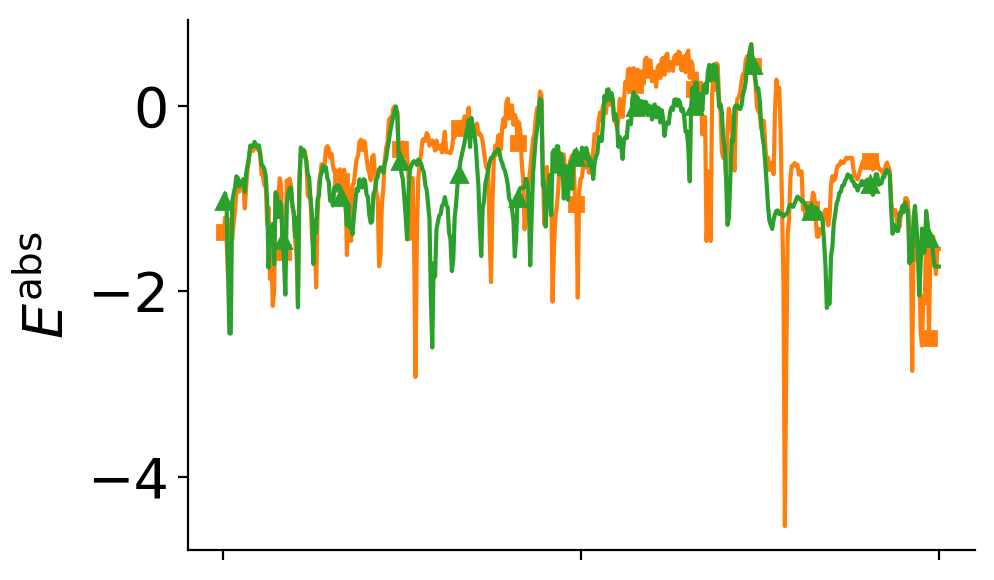}
    \end{subfigure}

    \caption{The background-error covariance is SOAR-correlated in the generation but misspecified as diagonal in the objective, in the cycled setting with $T=5$. (top) Truth profile and mean reconstructions of the representative variable obtained by Gaussian, the TV prior transform, and the residual prior transform for (left) Burgers, (middle left) Stoker, (middle right) Sod, and (right) Shu-Osher. (bottom) The corresponding pointwise error \eqref{eq:pw_error} on a $\log_{10}$ scale.}
    \label{fig:mismatch_cycled}
\end{figure}

Using the Stoker problem as a representative example, we also show the empirical $95\%$ band of the analysis over the 100 Monte Carlo realizations, given by the 2.5$th$ and 97.5$th$ percentiles. The band is computed pointwise and quantifies the spread of the analysis across the realizations. \Cref{fig:band_mismatch_stoker} shows the mean reconstruction together with the empirical $95\%$ band for the Stoker problem. The mean reconstructions of all three methods are comparable and stay close to the true profile, but the $95\%$ bands differ substantially. The unregularized method has the widest band, while both regularized methods in \eqref{eq:rvda_cost} are much tighter. Over the rarefaction wave in $h$, the residual prior transform gives the narrowest band, whereas the staircasing of the first-order local differencing produces a wider one.

\Cref{fig:mismatch_cycled} shows the results under the same misspecification but in the cycled setting with $T=5$ in \eqref{eq:forecast_model}. The unregularized reconstruction from \eqref{eq:3dvar_cost} exposes more of its instability as noise accumulates over the cycles. It survives only the Burgers problem, which has the simplest behavior of the four. Both regularized methods remain usable for all four problems. However, the TV prior transform yields strong oscillations and the solution departs significantly from the true profile for the Shu-Osher problem. In contrast, the residual prior transform solution stays smooth except over a small region in the Sod problem, and it captures the essential structure of the Shu-Osher problem, where the first-order local differencing fails badly over the region of rich variation.

\subsection{Sensitivity analysis}
\label{sec:sensitivity}

\begin{figure}[!ht]
    \centering
    \includegraphics[width=\textwidth]{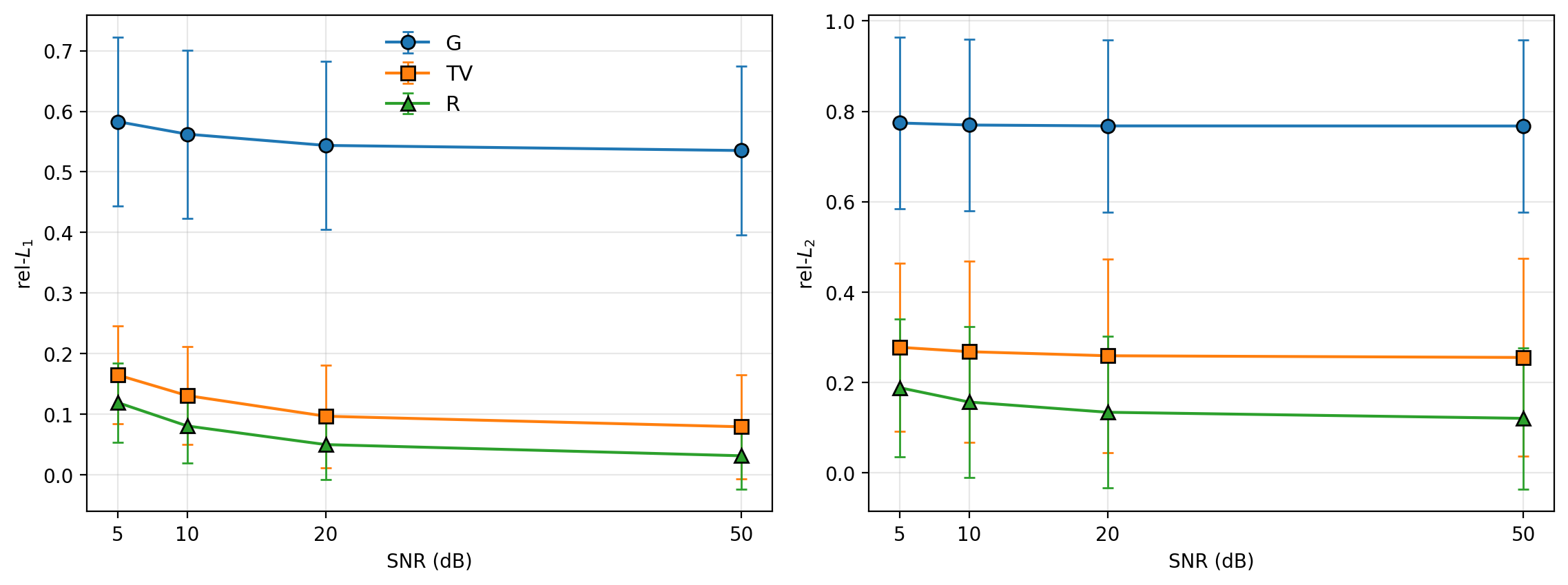}
    \caption{Sensitivity to the signal-to-noise ratio.}
    \label{fig:sweep_snr}
\end{figure}

Because Burgers' equation forms clear jump discontinuities at the shock location, we use it as the test bed for a sensitivity analysis of our approach under varying signal-to-noise ratio, grid resolution, and observation stride. We further include a sweep over covariance inflation, a classical remedy in both variational \cite{bonavita2012use,tabeart2020improving} and ensemble-based \cite{anderson1999monte, anderson2007adaptive,anderson2009spatially} methods, to verify that the improvement from the residual prior transform cannot be recovered by inflating the covariance alone. Multiplicative inflation adjusts only the amplitude and cannot alter the correlation structure \cite{tabeart2020improving}, so it cannot compensate for a structural misspecification.

\begin{figure}[!ht]
    \centering
    \includegraphics[width=\textwidth]{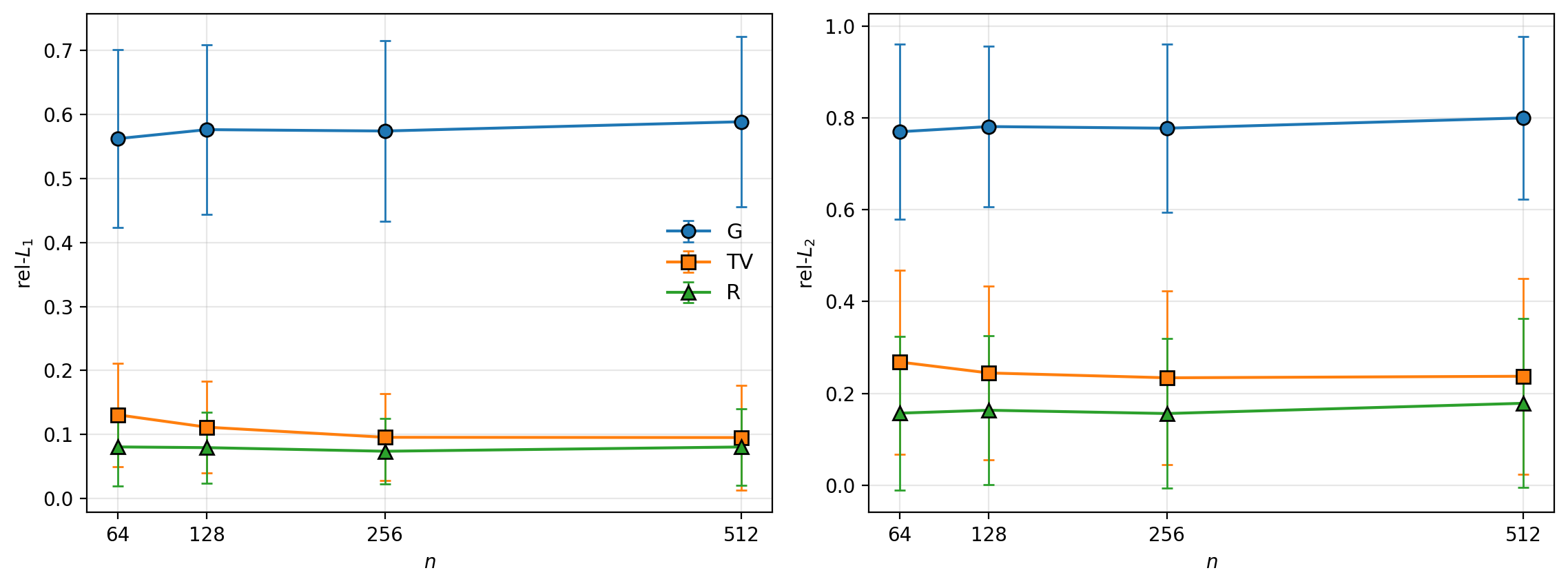}
    \caption{Sensitivity to the grid resolution.}
    \label{fig:sweep_nx}
\end{figure}

\begin{figure}[!ht]
    \centering
    \includegraphics[width=\textwidth]{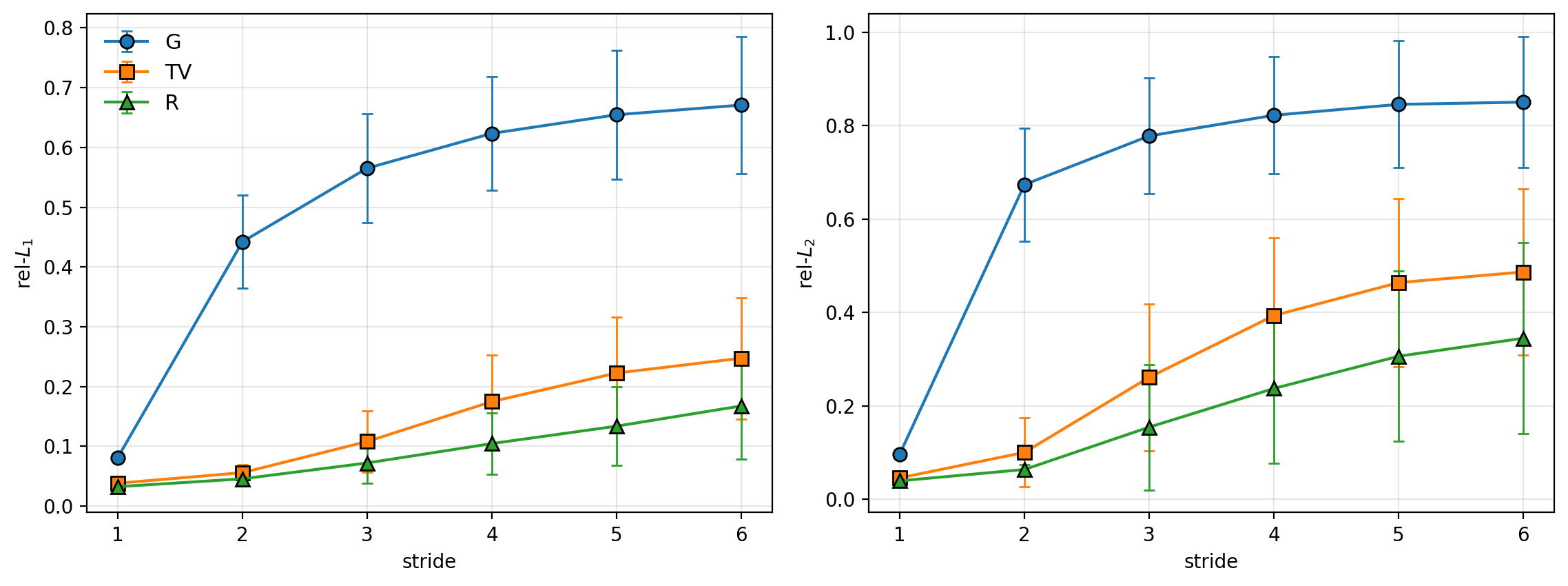}
    \caption{Sensitivity to observation sparsity.}
    \label{fig:sweep_stride}
\end{figure}

We evaluate each analysis against the truth using relative errors in both the $\ell_1$ and $\ell_2$ norms. For a single field, let $\hat{\bm f}$ denote the analysis and $\bm f$ the truth on the solution grid. The relative errors are
\begin{equation}\label{eq:rel_errors}
\text{rel-}L_1 = \frac{\sum_i \abs{\hat{f}_i - f_i}}{\sum_i \abs{f_i}}, \qquad
\text{rel-}L_2 = \sqrt{\frac{\sum_i (\hat{f}_i - f_i)^2}{\sum_i f_i^2}}.
\end{equation}
For a system with several variables, we compute the error of each variable separately rather than stacking them into a single vector. The relative $\ell_2$ error measures the overall performance of each method, while the relative $\ell_1$ error gives a closer look at the sharp and discontinuous features.

\begin{figure}[!ht]
    \centering
    \includegraphics[width=\textwidth]{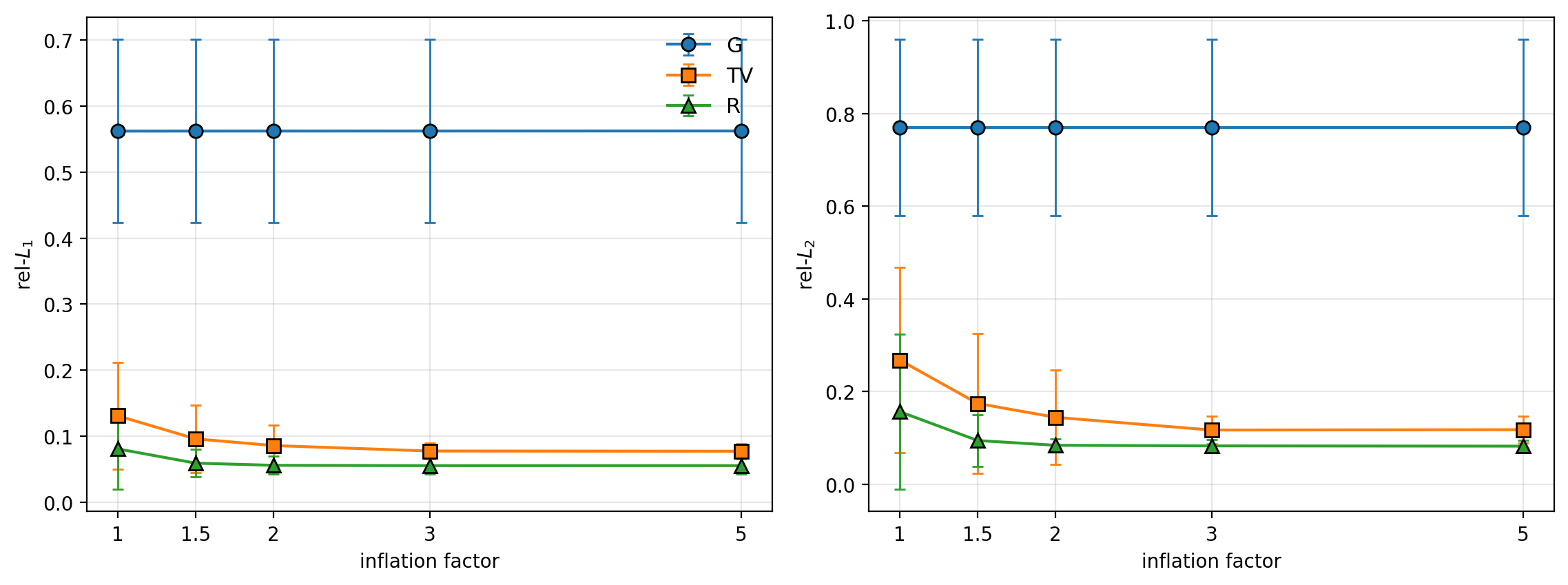}
    \caption{Effect of inflating the assumed diagonal background-error covariance.}
    \label{fig:sweep_inflation}
\end{figure}

All the sensitivity results are obtained from $100$ Monte Carlo simulations, using the misspecified background-error covariance in the single-cycle setting. For each realization we compute the relative error of that realization, and we report the mean and the $\pm1$ standard deviation of these per-realization errors as the curve and the error bar. This differs from the pointwise error in \eqref{eq:pw_error}, which is the error of the mean analysis, whereas here we average the error of each individual realization.

\Cref{fig:sweep_stride,fig:sweep_snr,fig:sweep_nx} report the relative errors as the observation stride, the signal-to-noise ratio, and the grid resolution are varied. In all three sweeps the two regularized methods are far more accurate than the Gaussian baseline, and the residual prior transform is the most accurate throughout. The relative errors are nearly flat under changes in the grid resolution, and they change only mildly with the signal-to-noise ratio. In these two sweeps the total variation and residual priors stay close to each other. In contrast, as the observations become sparser, the advantage of the residual prior transform over the total variation prior grows steadily, so the residual prior transform is especially effective for spatially sparse observations.

\Cref{fig:sweep_inflation} examines whether inflating the misspecified diagonal covariance in the objective can recover the accuracy lost to the structural misspecification. Inflating the diagonal covariance leaves the Gaussian baseline essentially unchanged, and although it improves the two regularized methods, their errors saturate well above the level attainable with a correctly specified covariance. Inflation adjusts only the amplitude of the covariance and cannot restore its correlation structure \cite{tabeart2020improving}. The improvement from the residual prior transform therefore cannot be reproduced by inflation alone.

\section{Conclusion}
\label{sec:conclusion}

This investigation develops a structure-preserving framework for the variational data assimilation of hyperbolic conservation laws with shocks, either forming from smooth initial data or present in discontinuous initial data. The framework integrates a residual prior transform into the regularized 3D-Var objective, solved with a hierarchical sparse Bayesian learning algorithm, which substantially improves on the unregularized $\ell_2$ analysis. Unlike the standard first-order local differencing prior, the residual prior transform does not commit to a fixed order of smoothness or assume such knowledge a priori. Numerical experiments on four test problems, under both correctly specified and misspecified background-error covariances, demonstrate that the residual prior transform is a more robust structure-preserving regularizer than first-order local differencing, and that its advantage grows as the assimilation conditions deteriorate. In both the correctly specified and misspecified cases, the first-order local differencing prior produces spurious oscillations, whereas the residual prior transform remains smoother, even when both priors stay comparably close to the true profile. The discrepancy between the two priors grows markedly over multiple assimilation cycles and as the observations become sparser.

Our framework opens several directions for future work. A natural extension is to two-dimensional conservation laws.  The residual transform $R = L_1 - L_2$ requires only that $L_1{\bm f} \approx L_2\bm f$ for locally smooth underlying function $f$, with disagreement signaling a nearby discontinuity. Because this criterion is a pointwise scalar comparison rather than a combination of oriented derivative components, it does not depend on the orientation of the underlying discontinuity relative to the grid. This avoids the isotropic/anisotropic combination question that arises in multidimensional TV, where the individual derivative components being combined are themselves orientation-dependent. On a rectangular grid, this permits a direct sequential extension of the residual transform by applying $L_1$ and $L_2$ along each grid direction, without the axis-aligned bias that motivates the isotropic/anisotropic distinction for TV.  A second direction concerns partial observations, where only a subset of the state variables is observed. Because the different state variables of the same dynamical system share a common sparsity profile, this joint sparsity could be leveraged to recover all state variables when only partial observations are available. The framework can also accommodate background-error covariance structures beyond the SOAR model used here, since the regularized objective does not depend on a particular choice. Finally, the approach extends naturally to the four-dimensional variational setting, in which observations distributed over a time window are assimilated together.

\bibliographystyle{siamplain}

\bibliography{references}
\end{document}